# ZERO TEMPERATURE LIMIT FOR THE BROWNIAN DIRECTED POLYMER AMONG POISSONIAN DISASTERS

BY RYOKI FUKUSHIMA[1] AND STEFAN JUNK[2]

*Kyoto University and Technical University of Munich*

We study a continuum model of directed polymer in random environment. The law of the polymer is defined as the Brownian motion conditioned to survive among space-time Poissonian disasters. This model is well studied in the positive temperature regime. However, at zero-temperature, even the existence of the free energy has not been proved. In this article, we show that the free energy exists and is continuous at zero-temperature.

**1. Introduction and main results.** We discuss the zero-temperature limit for a model of directed polymer in random environment. This work is partially motivated by the recent work [13] which studies the number of the open paths in the oriented percolation. In the directed polymer context, the main result in [13] corresponds to the existence of the free energy at zero-temperature for the Bernoulli environment. What makes this problem nontrivial is that at zero-temperature, the finite volume free energy has infinite mean, and hence the standard subadditivity argument fails. The proof in [13] instead relies on an intricate combination of various techniques developed in the study of the contact process and the oriented percolation. It would be desirable to know whether the zero-temperature model can be approximated by the more tractable low positive temperature model. However, to our knowledge, it is not known whether the free energy is continuous at zero-temperature.

In this article, we investigate a certain time-space continuous analogue of the model described above which has a similar nonintegrability issue at zero temperature. We prove the existence and the continuity of the free energy. The argument for the existence is a suitable modification of the standard subadditivity argument which is quite different from [13].

1.1. *The model.* We study the Brownian directed polymer in Poissonian environment introduced in [9]. We first recall the model at positive temperature and then introduce a natural zero-temperature version. Let us denote by $(B = \{B(t) : t \geq 0\}, P)$ the standard $d$-dimensional Brownian motion stating at the origin, and by $(\omega, \mathbb{P})$ the Poisson point process independent of $B$ with unit intensity

Received October 2018.
[1]Supported by JSPS KAKENHI Grant Number 16K05200.
[2]Supported by TUM Graduate School Internationalization Support.
*MSC2010 subject classifications.* Primary 60K37; secondary 60K35, 82A51, 82D30.

on $\mathbb{R}_+ \times \mathbb{R}^d$ ($d \geq 1$). The process $\omega$ is realized as a locally finite point measure as usual but with some abuse of notation, we will frequently identify $\omega$, and more generally any point measure, with its support. Let $U(x) \subset \mathbb{R}^d$ be the ball of unit volume centered at $x$ and $V_t$ denote a tube around the path of $B$:

$$V_t(B) := \{(s, x) \in [0, t) \times \mathbb{R}^d : x \in U(B(s))\}.$$

Then for given $\beta > 0$, the so-called polymer measure is defined as

(1.1) $$dP_t^{\beta,\omega} = \frac{1}{Z_t^{\beta,\omega}} \exp(-\beta \omega(V_t)) \, dP,$$

where

(1.2) $$Z_t^{\beta,\omega} = E[\exp(-\beta \omega(V_t))]$$

is the normalizing constant. Under this measure, the polymer receives a *repulsive* interaction from a point $(s, x) \in \omega$.

REMARK 1.1. In the earlier works [6, 8–11] on Brownian directed polymers, there is no negative sign in front of $\beta$ in the definitions (1.1) and (1.2) and general $\beta \in \mathbb{R}$ is considered. We use the above formulation since we focus on the repulsive case in this article.

In our repulsive case, the above polymer measure can be interpreted as the law of the Brownian motion conditioned to survive among Poissonian disasters. More precisely, we enlarge the probability space for the Brownian motion and introduce an Exp(1) random variable $\xi$, independent of $B$ and $\omega$, and define the "death time"

(1.3) $$\tau_\beta(\omega) := \inf\{t \geq 0 : \beta \omega(V_t) \geq \xi\}.$$

Then we have

(1.4) $$Z_t^{\beta,\omega} = P(\tau_\beta(\omega) \geq t)$$

and thus the polymer measure admits the aforementioned interpretation. Now in the zero-temperature limit $\beta \to \infty$, the above death time becomes

(1.5) $$\tau_\infty(\omega) := \inf\{t \geq 0 : \omega(V_t) \geq 1\},$$

that is, the process is immediately killed when it gets close to a disaster. Note that this stopping time can alternatively expressed as the hitting time to the set of disasters defined by

(1.6) $$\mathcal{D} := \bigcup_{(s,x) \in \omega} \{s\} \times U(x).$$

1.2. *Free energy.* As is usual in the study of a model of statistical physics, it is important to understand the asymptotics of the normalizing constant $Z_t^{\beta,\omega}$. For positive temperature, the following result is known.

THEOREM A (Theorem 2.2.1 in [9]). *There exists a continuous function $p\colon [0,\infty) \to (-\infty, 0]$ such that for almost all $\omega$,*

$$(1.7) \quad \lim_{t\to\infty} \frac{1}{t} \log P(\tau_\beta(\omega) \geq t) = \lim_{t\to\infty} \frac{1}{t} \mathbb{E}\big[\log P(\tau_\beta(\omega) \geq t)\big] = p(\beta).$$

The limit $p(\beta)$ is an important observable called the free energy and it is for instance known to characterize a localization transition of the polymer; see [9, 10] for more detail.

The goal of this article is to extend the above existence and continuity to the value $\beta = \infty$. The methods of proving the above theorem do not seem to cover the case $\beta = \infty$: First, the existence of the limit is proved by the superadditivity of the mean $\mathbb{E}[\log P(\tau_\beta(\omega) \geq t)]$ and a concentration bound for $\log P(\tau_\beta(\omega) \geq t)$ around its mean. However, as we will see in Proposition 1.2, the mean $\mathbb{E}[\log P(\tau_\infty(\omega) \geq t)]$ does not exist. One could alternatively use the subadditive ergodic theorem but this approach also requires the integrability. Second, the continuity for $\beta \in (0,\infty)$ is a consequence of the convexity of $p(\cdot)$, which essentially follows from the Hölder inequality. But the convexity tells us nothing about the continuity at the boundary $\beta = \infty$. In order to motivate the notation used in the statement of the main result, let us observe why the integrability of $\log P(\tau_\infty(\omega) \geq t)$ is violated.

PROPOSITION 1.2. *For any $t > 0$, $\mathbb{E}[\log P(\tau_\infty(\omega) \geq t)] = -\infty$.*

PROOF. Let $F$ be the first disaster time close to the origin:

$$F := \inf\left\{t \geq 0 : \omega\left([0,t] \times \frac{1}{2}U(0)\right) \neq 0\right\}.$$

Note that the Brownian motion gets killed if $B(F) \in \frac{1}{2}U(0)$. Thus on $\{F < t\}$, we have

$$P(\tau_\infty(\omega) \geq t) \leq P\left(B(F) \notin \frac{1}{2}U(0)\right) \leq \exp\left(-\frac{C}{F}\right).$$

Since $F$ is exponentially distributed, we have

$$\mathbb{E}\big[\log P(\tau_\infty(\omega) \geq t)\big] \leq -C\mathbb{E}\big[F^{-1}\mathbb{1}\{F < t\}\big] = -\infty. \qquad \square$$

The proof of this proposition suggests that the nonintegrability is caused by the existence of a Poisson point near the starting point of the Brownian motion. It is reasonable to believe that this is the *only* source of nonintegrability and we will in fact confirm this intuition in the proof.

1.3. *Main results.*  In view of the discussion at the end of the previous subsection, it is natural to consider a modified death time where we ignore the disasters in the first unit time interval. To this end, for $I \subset \mathbb{R}_+$, we write $\omega_I$ for the restriction $\omega|_{I \times \mathbb{R}^d}$ as a measure and define

$$\tau_\beta^1(\omega) := \tau_\beta(\omega_{[0,1]^c}) = \begin{cases} \inf\{t \geq 1 \colon \beta\omega_{[0,1]^c}(V_t) \geq \xi\} & \text{for } \beta < \infty, \\ \inf\{t \geq 1 \colon \omega_{[0,1]^c}(V_t) \geq 1\} & \text{for } \beta = \infty. \end{cases}$$

It is often convenient to restrict the Brownian motion to a domain growing at polynomial speed:

(1.8) $$\mathcal{A}_t := \left\{ \sup_{0 \leq s \leq t} |B(s) - B(0)| \leq \lceil t \rceil^2 \right\}.$$

The probability of $\mathcal{A}_t^c$ is bounded by $\exp(-ct^3)$ by the reflection principle, and hence it should be much smaller than the survival probability.

Now we are ready to state the main theorem.

THEOREM 1.3.  *There exists $p(\infty) \in (-\infty, 0]$ such that the following hold*:
  (i) $\lim_{t \to \infty} \frac{1}{t} \mathbb{E}[\log P(\tau_\infty^1(\omega) \geq t, \mathcal{A}_t)] = p(\infty)$,
  (ii) *for almost all $\omega$,* $\lim_{t \to \infty} \frac{1}{t} \log P(\tau_\infty(\omega) \geq t) = p(\infty)$,
  (iii) $\lim_{\beta \to \infty} p(\beta) = p(\infty)$.

REMARK 1.4.  The proof of this theorem is almost identical for $d = 1$ and $d \geq 2$, except for one point which we mention in Remark 2.6. For this reason, we carry out the proof mostly in the one-dimensional setting and then describe the necessary modification to deal with higher dimensions in the final section.

1.4. *Related works.*  For the general background and known results on the directed polymer in random environments, we refer the reader to [5].

The Brownian directed polymer in Poissonian environment was first introduced and studied in [9]. One of the advantages of this model is that one can employ tools from stochastic analysis. In addition to standard results, such as the existence and phase transition of the free energy, some results on the fluctuation of the free energy and the displacement of the polymer were proved. The latter results are sharper than what is known for the models based on the simple random walk. Later in [10], a more precise localization result for the polymer was established. Our continuity result for the free energy suggests—but does not prove—that these results remain true at zero-temperature. This direction is worthy of further investigation.

Next, we turn to other works on the zero-temperature limit. There are not so many results on the zero-temperature limit for the directed polymer in random environment. This is mainly because we have a simple answer in a large class of settings. To see this, let us consider the most studied nearest neighbor lattice polymer setting. In this case, the environment is given by real-valued random variables $\omega = (\omega(k,x))_{k\in\mathbb{N}, x\in\mathbb{Z}^d}$ and the polymer measure is defined as follows:

$$(1.9) \qquad P_n^{\beta,\omega}(\gamma) = \frac{1}{Z_n^{\beta,\omega}} \exp\left(-\beta \sum_{k=1}^n \omega(k,\gamma(k))\right) \mathbb{1}\{\gamma \in \mathcal{N}_n\},$$

where $\mathcal{N}_n$ denotes the set of nearest neighbor paths of length $n$ on $\mathbb{Z}^d$. Now if the *time constant* for the directed first passage percolation

$$(1.10) \qquad \mu = \lim_{n\to\infty} \frac{1}{n} \min_\gamma \sum_{k=1}^n \omega(k,\gamma(k))$$

is nonzero (this holds, e.g., when $\operatorname{ess\,inf}\omega \neq 0$), then it is easy to deduce a continuity result

$$(1.11) \qquad \lim_{\beta\to\infty}\lim_{n\to\infty} \frac{1}{\beta n} \log Z_n^{\beta,\omega} = -\mu.$$

On the other hand, if $\operatorname{ess\,inf}\omega = 0$ and the set $\{(k,x): \omega(k,x) = 0\}$ percolates, then we have $\mu = 0$. In this case, $Z_n^{\infty,\omega} = \lim_{\beta\to\infty} Z_n^{\beta,\omega}$ is the number of open paths and $\lim_{n\to\infty} \frac{1}{n} \log Z_n^{\infty,\omega}$ represents the growth rate. As is mentioned at the beginning of this article, the existence of this limit is proved in [13] but it is not known whether the limit equals $\lim_{\beta\to\infty} \lim_{n\to\infty} \frac{1}{n} \log Z_n^{\beta,\omega}$. Two recent works [7, 15] study this type of problem in a nonnearest neighbor model on $\mathbb{Z}_+ \times \mathbb{Z}^d$ defined by

$$(1.12) \qquad P_n^{\beta,\omega}(\gamma) = \frac{1}{Z_n^{\beta,\omega}} \exp\left(-\sum_{k=1}^n [\beta\omega(k,\gamma(k)) + |\gamma(k-1) - \gamma(k)|^\alpha]\right)$$

and proved the continuity of the free energy at $\beta = \infty$. In this case, $\log Z_n^{\beta,\omega}$ is integrable, and hence the existence follows from the subadditivity argument.

Finally, there is a recent work [2] where the zero-temperature limit of the polymer measure is discussed for the model on $\mathbb{Z}_+ \times \mathbb{R}$ defined by

$$(1.13) \qquad \begin{aligned} & P_n^{\beta,\omega}(d\gamma) \\ & = \frac{1}{Z_n^{\beta,\omega}} \exp\left(-\beta \sum_{k=1}^n [\omega(k,\gamma(k)) + |\gamma(k-1) - \gamma(k)|^2]\right) \prod_{k=1}^n d\gamma_k. \end{aligned}$$

In the preceding works [1] and [3], the infinite volume polymer measure is constructed for every given asymptotic slope, at zero and positive temperature, respectively. Then in [2], it is proved that as $\beta \to \infty$, not only the free energy but also the infinite volume polymer measure converges. This model has a similarity to the model studied in this article since the polymer measure in (1.1) has a heuristic representation

$$(1.14) \quad P_t^{\beta,\omega}(\mathrm{d}\gamma) = \frac{1}{Z_t^{\beta,\omega}} \exp\left(-\beta\omega(V_t) - \frac{1}{2}\int_0^t |\dot{\gamma}(s)|^2 \mathrm{d}s\right) \mathrm{d}\gamma.$$

However, we do not multiply the term $\int_0^t |\dot{\gamma}(s)|^2 \mathrm{d}s$ by $\beta$, and thus the two models behave quite differently as $\beta \to \infty$. The zero temperature model in [1] is of last passage percolation type and concentrates on a single path, whereas our result implies that the entropy is nondegenerate at zero temperature.

**2. High-level structure of proof.** In this section, we explain the high-level structure of the proof. In order to make the flow of the argument clear, we discuss the convergence only for $t \in \mathbb{N}$. The complete proof will be given in Section 6.

*Convergence of the mean.* We first need to prove that $\mathbb{E}[\log P_\omega(\tau_\infty^1 \geq t)]$ is finite. This will follow as a corollary to Lemma 3.1 in Section 3—we refrain from stating it precisely since it is designed to cover a more complicated situation and requires a number of terminologies. Now if in addition $\{\mathbb{E}[\log P(\tau_\infty^1(\omega) \geq t)]\}_{t\geq 0}$ were a superadditive sequence, then Theorem 1.3(i) would follow. However, the modification $\tau_\beta \to \tau_\beta^1$ makes it difficult to prove the superadditivity. The standard argument for the superadditivity (see, e.g., [9], Section 6) yields that, for $s > 1$,

$$(2.1) \quad \begin{aligned} &\mathbb{E}[\log P(\tau_\infty^1(\omega) \geq s+t)] \\ &\geq \mathbb{E}[\log P(\tau_\infty^1(\omega) \geq s)] + \mathbb{E}[\log P(\tau_\infty(\omega) \geq t)], \end{aligned}$$

in which we get $\tau_\infty$ instead of $\tau_\infty^1$ in the second term. For this reason, we shall instead prove the following *almost* superadditivity that is known to imply Theorem 1.3(i) by [14], Theorem 2.

PROPOSITION 2.1. *Let $a_\beta(t) := \mathbb{E}[\log P(\tau_\beta^1(\omega) \geq t, \mathcal{A}_t)]$. For every $\delta \in (0, \frac{1}{2})$, there exists $t_0 > 0$ independent of $\beta$ such that for all $s, t \geq t_0$,*

$$(2.2) \quad a_\beta(s+t) \geq a_\beta(s) + a_\beta(t) - (s+t)^\delta.$$

*Almost sure convergence with the modification.* In order to upgrade the convergence of the mean to the almost sure convergence, we prove the following concentration bound.

PROPOSITION 2.2. *For every $\delta \in (0, \frac{1}{2})$ and all $r \geq 0$, there exists $t_0 > 0$ such that for all $t \geq t_0$ and $\beta \in [0, \infty]$,*

$$(2.3) \quad \mathbb{P}(|\log P(\tau_\beta^1(\omega) \geq t, \mathcal{A}_t) - \mathbb{E}[\log P(\tau_\beta^1(\omega) \geq t, \mathcal{A}_t)]| \geq t^{\frac{1}{2}+\delta}) \leq t^{-r}.$$

REMARK 2.3. For fixed $\beta < \infty$, an exponential concentration bound is obtained in Theorem 2.4.1(b) in [9]. However, it does not cover the case $\beta = \infty$ since it contains a constant that degenerates at $\beta = \infty$.

Proposition 2.2 together with Theorem 1.3(i) implies the almost sure convergence of $t^{-1} \log P(\tau_\beta^1(\omega) \geq t, \mathcal{A}_t)$.

*Continuity of the free energy.* Given the concentration bound in Proposition 2.2, we can derive the following estimate on the rate of convergence by adapting the argument in [17] for first passage percolation.

PROPOSITION 2.4. *For every $\delta > 0$, there exists $t_0 > 0$ such that for all $t \geq t_0$ and $\beta \in [0, \infty]$,*

$$\left|\mathbb{E}[\log P(\tau_\beta^1(\omega) \geq t, \mathcal{A}_t)] - tp(\beta)\right| \leq t^{\frac{1}{2}+\delta}. \tag{2.4}$$

We can therefore approximate $p(\beta)$ by $\frac{1}{t}\mathbb{E}[\log P(\tau_\beta^1(\omega) \geq t, \mathcal{A}_t)]$ with large $t$, uniformly over all $\beta \in [0, \infty]$, and the continuity from Theorem 1.3(iii) follows because this expectation depends only on the disasters in a bounded subset of $\mathbb{R}_+ \times \mathbb{R}$.

*Getting rid of the modification.* It remains to remove the modification in the time interval $[0, 1]$. This might look an easy task but in fact is a subtle problem. It is possible to replace $[0, 1]$ by $[0, \epsilon]$ for any $\epsilon > 0$ and prove that the limit of $\frac{1}{t}\mathbb{E}[\log P(\tau_\beta^\epsilon(\omega) \geq t, \mathcal{A}_t)]$ is independent of $\epsilon$. If we knew in addition that

$$P_t^{\beta,\omega}(B(u) \in (-R, R) \text{ for all } u \in [0, 1]) = e^{o(t)} \tag{2.5}$$

for large $R > 0$, then we could restrict the consideration to the above event and argue that there are no disasters in $[0, \epsilon] \times [-R, R]$ for sufficiently small $\epsilon > 0$, which implies $\tau_\beta^\epsilon(\omega) = \tau_\beta(\omega)$. However, proving (2.5) turns out to be difficult. We instead prove that the disasters in the time interval $[0, 1]$ do not affect the survival probability too much, uniformly in the endpoint.

PROPOSITION 2.5. *There exists a finite positive random variable $A(\omega)$ such that for all $x \in \mathbb{R}$,*

$$P(B(2) \in dx, \tau_\infty(\omega) \geq 2) \geq A(\omega) P(B(1) \in dx).$$

REMARK 2.6. This is the point where we need an extra argument for the higher dimensions. In fact, the above proposition fails to hold as it is for $d \geq 3$. We defer the description of the extra argument to the end of the article since it requires several auxiliary definitions. We include the case $d = 2$ there in order to make the notation simple in the other sections.

*Key technical steps.*   The main technical difficulty lies in the proofs of Propositions 2.1 and 2.2, which share much in common. Indeed, the former consists of controlling the effect of removing disasters in $[t, t+1] \times \mathbb{R}$ to the survival probability, and the latter relies on the control on the *influence* of resampling the disasters in $[i, i+1] \times \mathbb{R}$, as is usual for concentration bounds. The following proposition provides those controls.

PROPOSITION 2.7.   *For every $p \in \mathbb{N}$, there exists $C > 0$ such that for all $\beta \in [0, \infty]$, all $t \geq C$ and $r, s > 0$ such that $1 \leq r \leq r + s \leq t$ and either $r + s \leq t - 1$ or $r + s = t$,*

$$
\mathbb{E}[|\log P(\tau_\beta^1(\omega) \geq t \mid \tau_\beta^1(\omega_{[r,r+s]^c}) \geq t, \mathcal{A}_t)|^p]
$$
(2.6)
$$
= \mathbb{E}[|\log P(\tau_\beta^1(\omega) \geq t, \mathcal{A}_t) - \log P(\tau_\beta^1(\omega_{[r,r+s]^c}) \geq t, \mathcal{A}_t)|^p]
$$
$$
\leq C(1 + s^p) + C(1 + \log^+ t)^C.
$$

The proof of Proposition 2.7 is the technical core of this work and will take up a large portion of the rest of the paper. We prove it by sampling the paths according to the conditional law $P(\cdot \mid \tau_\beta^1(\omega_{[r,r+s]^c}) \geq t, \mathcal{A}_t)$ and then estimating the cost for the paths to avoid the additional disasters in $[r, r+s] \times \mathbb{R}$. This requires a lower bound on the survival probability for the Brownian motion whose initial and terminal distribution at time $r$ and $r + s$ are given by the above conditional law.

We will prove such a lower bound in three steps. First, in Lemma 3.1, we prove a lower bound on the survival probability for the Brownian bridge with a further additional restriction that it stays in a tube around the line connecting its initial and terminal points (Section 3). Then in Lemma 4.1, we prove a moment bound for the general initial and terminal distribution (Section 4). This is done by *duplicating* the tube strategy provided by Lemma 3.1. In order to make use of many tubes, we need the endpoint distribution to be dispersed (see Figure 1). Finally, in Lemma 5.1, we show that the distribution of $(B(r), B(r+s))$ under the conditional law $P(\cdot \mid \tau_\beta^1(\omega_{[r,r+s]^c}) \geq t, \mathcal{A}_t)$ is indeed dispersed (Section 5).

Sections 3, 4 and 5 are isolated in the sense that nothing beyond the main lemma is used later on.

*Notational convention.*   In the proof, we use $c$ and $C$ to denote positive constants whose values may change from line to line.

**3. Survival probability in a tube.**   In this section, we provide a lower bound for the survival probability of the Brownian motion which is conditioned to end at a fixed point and restricted to a tube.

We start by introducing some notation. To describe the first conditioning, we write $P^{r,x;s,y}$ for the Brownian bridge measure between $(r, x)$ and $(s, y)$. We introduce the following intervals:

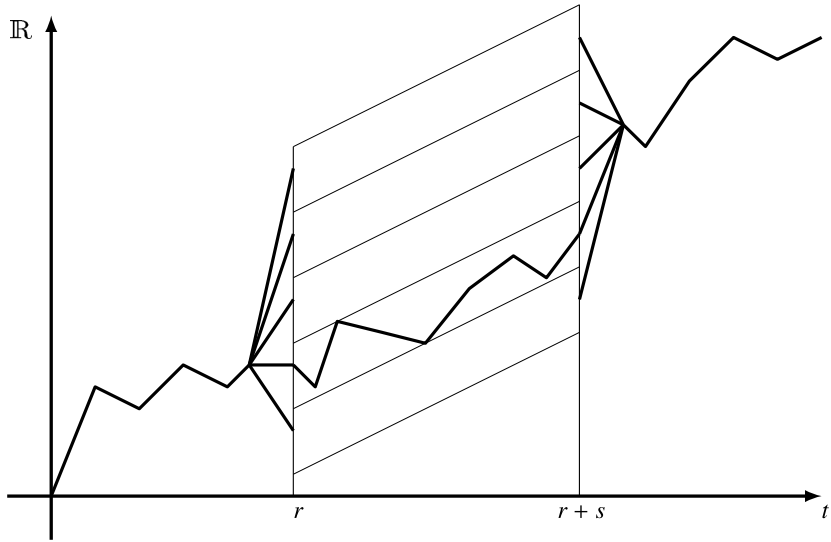

FIG. 1. *An illustration of the proof strategy of Proposition* 2.7. *The survival probability in each tube is controlled by Lemma* 3.1. *If the distribution of the polymer (conditioned to avoid the disasters outside* $[r, r+s]$) *at times* $r$ *and* $r+s$ *is sufficiently dispersed, then we can choose the best tube among many in order to construct a good survival strategy on the whole time interval* $[0, t]$.

- $x$ and $y$ will be chosen from $J^{(5)} := [-\frac{5}{2}, \frac{5}{2}]$,
- the Brownian motion will be restricted to $J^{(6)} := [-3, 3]$,
- then the survival probability depends only on the disasters in $J^{(7)} := [-\frac{7}{2}, \frac{7}{2}]$,
- $J^{(5)}$ is divided into $\mathcal{J} = \{J^{(1)}_{-2}, \ldots, J^{(1)}_{2}\}$, where $J^{(1)}_x := x + [-\frac{1}{2}, \frac{1}{2})$.

The role of $J^{(7)}$ is to ensure the independence of the survival probabilities in different tubes in our duplication strategy in Figure 1. For $t > 0$, let $F_t$ denote the first disaster in $[0, t] \times J^{(7)}$, that is,

$$F_t := \inf\{r \in [0, t] : \exists z \in J^{(7)} \text{ such that } (r, z) \in \omega\}$$

with the convention $F_t = t$ if $\omega \cap [0, t] \times J^{(7)} = \emptyset$. Similarly, we let

$$L_t := \sup\{r \in [0, t] : \exists z \in J^{(7)} \text{ such that } (r, z) \in \omega\}$$

denote the last disaster in $[0, t] \times J^{(7)}$, where we set $L_t = 0$ if there is no such disaster. The goal of this section is the following lemma, which provides a lower bound on the survival probability in the tube $[0, t] \times J^{(6)}$:

LEMMA 3.1. *There exists* $C > 0$ *such that the following hold* $\mathbb{P}$-*almost surely*:

(i) *For all* $x, y \in J^{(5)}$,

(3.1)
$$\mathbb{E}[\log P^{0,x;t,y}(\tau_\infty(\omega) \geq t, B(s) \in J^{(6)} \text{ for all } 0 \leq s \leq t) \mid F_t, L_t]$$
$$\geq -C(t + \mathbb{1}\{F_t < t\}(F_t^{-1} + (t - L_t)^{-1})),$$

(ii) $\mathbb{E}[\log P(\tau^1_\infty(\omega) \geq t, B(s) \in J^{(6)} \text{ for all } 0 \leq s \leq t)] \geq -C(t + 1).$

REMARK 3.2. The tube is assumed to be parallel to the time axis but this is not restrictive, as we can change the terminal point of the Brownian bridge by applying a time-space affine transformation which leaves the law of $\omega$ invariant. We include this generalization to Lemma 4.1 since the proof of Lemma 3.1 in the above simple form is already quite long and complicated.

The terms $F_t^{-1}$ and $(t - L_t)^{-1}$ above are the costs for the Brownian motion to avoid the first and last disasters in $[0, t] \times J^{(7)}$, respectively. Therefore, this lemma justifies the intuition discussed after Proposition 1.2. To see the reason why the cost is inverse proportion of $F_t$ or $(t - L_t)$, we state simple estimates for Brownian motion without proof, which we will repeatedly use in the proof.

LEMMA 3.3. *There exists $C > 0$ such that for every $s, t > 0$ and $x, y \in \{-2, \ldots, 2\}$, almost surely on $\{B(t) \in J_x^{(1)}\}$,*

$$P\bigl(B(s + t) \in J_y^{(1)} \text{ and } B(u + t) \in J^{(6)} \text{ for all } u \in [0, s] \mid B(t)\bigr)$$

(3.2)
$$\geq \begin{cases} e^{-\frac{C}{s} - Cs} & \text{if } x \neq y, \\ e^{-Cs} & \text{if } x = y. \end{cases}$$

We are going to bound the probability in (3.1) from below by constructing a specific survival strategy for the Brownian motion. We will introduce various terminologies in the course of describing the strategy. Given an environment $\omega$, we can find $T_i \geq 0$ and $D_i \in J^{(7)}$ such that

$$\omega \cap (\mathbb{R}_+ \times J^{(7)}) = \{(T_0, D_0), (T_1, D_1), \ldots\}$$

and such that $T_0 < T_1 < \cdots$. We denote the interarrival times by $\Delta_0 := T_0$ and

$$\Delta_i := T_i - T_{i-1}$$

for $i \geq 1$, which are independent exponential random variables with parameter 7. We say that $J_x^{(1)} \in \mathcal{J}$ is *contaminated* by $(T_j, D_j)$ if

$$J_x^{(1)} \cap U(D_j) \neq \varnothing.$$

It is simple to check that if $J_x^{(1)} \in \mathcal{J}$ is not contaminated by $(T_j, D_j)$ and $B(T_j) \in J_x^{(1)}$, then the Brownian motion is not affected by the disaster at time $T_j$. Clearly, every disaster can contaminate at most two sites, and since $|\mathcal{J}| = 5$, there exists a sequence $(s(0), s(1), \ldots) \in \{0, 1, \ldots, 4\}^{\mathbb{N}}$ such that $J_{s(j)}^{(1)}$ is not contaminated by $(T_j, D_j)$ or $(T_{j+1}, D_{j+1})$; see Figure 2. The interval $J_{s(j)}^{(1)}$ is *safe* in the sense that the Brownian motion can survive during $[T_j, T_{j+2})$ simply by staying there.

Note that if there is no disaster in $[0, t] \times J^{(7)}$ (that is, on $\{F_t = t\} = \{F_t = t, L_t = 0\}$), we get (3.1) from Lemma 3.3 since

$$P(\tau(\omega) \geq t) \geq P\bigl(B(s) \in J^{(6)} \text{ for all } s \in [0, t]\bigr) \geq e^{-Ct}.$$

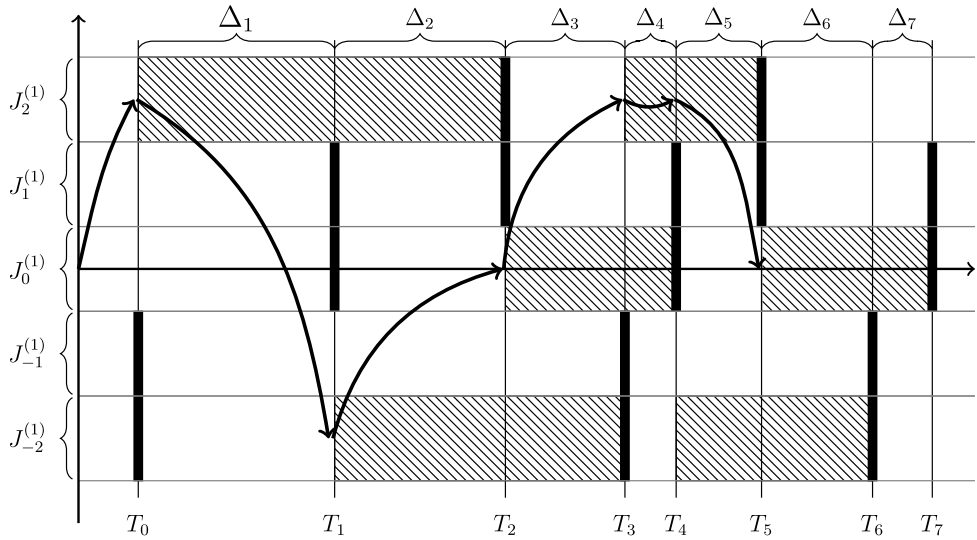

FIG. 2. *An illustration of the survival strategy until the first regeneration time $R_1$. In this figure, we have $\rho_1 = 5$. At every disaster time, (typically) two intervals are contaminated (marked by the thick lines). The left ends of the striped regions are safe intervals. The arrows indicate to which interval the Brownian motion is supposed to move.*

For the remainder of this section, we only discuss the case $\{F_t < t\} = \{F_t < t, L_t > 0\}$.

*The first interval.* The survival strategy up to $T_0 = F_t$ is prescribed by the event

(3.3) $$\mathcal{S}(0) := \{B(T_0) \in J^{(1)}_{s(0)} \text{ and } B(u) \in J^{(6)} \text{ for } u \in [0, T_0]\}.$$

From the estimates in Lemma 3.3, we get

$$\log P(\mathcal{S}(0)) \geq -C(F_t + F_t^{-1}).$$

*Renewal construction.* After $T_0 = F_t$, we define the sequence of survival strategies by using a renewal structure. Let $\rho_0 := 0$ and for $i \geq 0$,

$$\rho_{i+1} = \inf\{j > \rho_i + 1 : \Delta_j > \Delta_{j-1}\}.$$

We write the corresponding disaster time by

$$R_i := T_{\rho_i}.$$

We now recursively define events $\mathcal{S}(i)$ ($i \geq 1$) as follows: $B(u) \in J^{(6)}$ for all $u \in [R_{i-1}, R_i)$ and in addition,

(S1) $\quad B(T_j) \in J^{(1)}_{s(j)} \quad$ for $j = \rho_{i-1}, \ldots, \rho_i - 2$;

(S2) $\quad B(u) \in J^{(1)}_{s(\rho_i - 2)} \quad$ for $u \in [T_{\rho_i - 2}, T_{\rho_i - 1}]$;

(S3) $\quad B(T_{\rho_i}) \in J^{(1)}_{s(\rho_i)}.$

In words, the Brownian motion moves to the safe interval in each time interval $(T_j, T_{j+1})$ except for $j = \rho_i - 2$. Note that we may have $\rho_i = \rho_{i-1} + 2$ and then the step (S1) is to be skipped. The second step (S2) is possible in this case since we have $B(T_{\rho_i-2}) = B(T_{\rho_{i-1}}) \in J^{(1)}_{s(\rho_{i-1})}$ by the definition of $\mathcal{S}(0)$ ($i = 1$) and $\mathcal{S}(i-1)$ ($i \geq 2$). Now on the event $\{\rho_1 = k\}$ ($k \geq 2$), Lemma 3.3 yields

$$(3.4) \qquad \log P(\mathcal{S}(1) \mid \mathcal{S}(0)) \geq -C \sum_{\substack{i=1,\ldots,k \\ i \neq k-1}} \Delta_i^{-1} - C \sum_{i=1}^{k} \Delta_i.$$

It is important that the term $\Delta_{k-1}^{-1} = \max\{\Delta_1^{-1}, \ldots, \Delta_k^{-1}\}$ is omitted from the first sum on the right-hand side, due to the unusual strategy in (S2) above. Indeed, if that sum was taken over $1 \leq i \leq k$, it would be the sum of inverse exponential random variables, which is not $\mathbb{P}$-integrable. On the other hand, the other terms $\{\Delta_1^{-1}, \ldots, \Delta_{k-2}^{-1}, \Delta_k^{-1}\}$ gain one extra degree of integrability from the knowledge that they are the $k-1$ smallest members from the collection $\{\Delta_1^{-1}, \ldots, \Delta_k^{-1}\}$.

*Last interval.* It remains to prescribe the behavior after the last renewal time before time $t$. Let us denote by

$$N(s) := \sum_{i=1}^{\infty} \mathbb{1}\{T_i \leq s\} \quad \text{and}$$

$$M(s) := \sum_{i=1}^{\infty} \mathbb{1}\{R_i \leq s\}$$

the numbers of disasters and renewals up to time $s$, respectively. We further set

$\sigma := N(L_t) - M(L_t) =$ the number of disasters in $[R_{M(L_t)}, L_t] \times J^{(7)}$,

$U := L_t - R_{M(L_t)} =$ the duration from the last renewal to $L_t$.

Then the survival strategy in $[R_{M(L_t)}, t]$ is prescribed by the event $\mathcal{T}$ defined as follows: $B(u) \in J^{(6)}$ for all $u \in [R_{M(L_t)}, t]$ and in addition,

(S4) $\qquad B(T_j) \in J^{(1)}_{s(j)} \quad \text{for } j = M(L_t), \ldots, N(L_t) - 1,$

(S5) $\qquad B(u) \in J^{(1)}_{s(N_t-1)} \quad \text{for } u \in [T_{N(L_t)-1}, L_t)),$

(S6) $\qquad B(t) = y.$

In the case where the last disaster time $L_t$ is a renewal time, both (S4) and (S5) are to be skipped. In words, the strategy $\mathcal{T}$ for the terminal part is the same as for the previous cases except that we choose to remain in $J^{(1)}_{s(N(L_t)-1)}$ after the last disaster before $L_t$, regardless of whether a renewal occurs after $L_t$ or not. Then exactly as in (3.4), on the event $\{\sigma = n\}$, we have

$$\log P^{0,x;t,y}(\mathcal{T} \mid \mathcal{S}(0), \ldots, \mathcal{S}(M(L_t)))$$
$$\geq -C\left(\sum_{i=1}^{n-1} \Delta_i^{-1} + \sum_{i=1}^{n} \Delta_i + (t - L_t) + (t - L_t)^{-1}\right),$$

where the last term $(t - L_t)^{-1}$ appears since the Brownian motion has to move from $J^{(1)}_{s(N(L_t)-1)}$ to the endpoint $y$ during $[L_t, t]$. Note that since there is no renewal in $[R_{M(L_t)}, L_t]$, the strategy $\mathcal{T}$ makes the Brownian motion survive without moving in the shortest interval among $\{[T_j, T_{j+1}]\}_{j=M(L_t)}^{N(L_t)-1}$. Therefore, for the same reason as before, we can expect that the sum $\sum_{i=1}^{n-1} \Delta_i^{-1}$ gains an extra degree of integrability.

Collecting the above strategies, we define

$$\mathcal{S}_t := \mathcal{S}(0) \cap \bigcap_{i=1}^{M(L_t)} \mathcal{S}(i) \cap \mathcal{T}.$$

Then the probability that the Brownian motion survives in the tube $[0, t] \times J^{(6)}$ is bounded from below by

$$\log P^{0,x;t,y}(\tau_\infty(\omega) \geq t, B(s) \in J^{(6)} \text{ for all } 0 \leq s \leq t)$$
$$\geq \log P(\mathcal{S}_t)$$
(3.5)
$$= \log P(\mathcal{S}(0)) + \sum_{i=1}^{M(L_t)} \log P(\mathcal{S}(i) \mid \mathcal{S}(i-1))$$
$$+ \log P(\mathcal{T} \mid \mathcal{S}(0), \ldots, \mathcal{S}(M(L_t))).$$

*Expectation conditioned on* $\{R_i\}_{i \geq 1}$. We are going to bound the $\mathbb{P}$-expectation of the last line in (3.5) conditioned on $F_t$ and $L_t$. What makes the argument complicated is that the number of summands $M(L_t)$ is random, depending on $\{R_i\}_{i \geq 1}$. Thus we need to estimate $\mathbb{E}[\log P(\mathcal{S}(i) \mid \mathcal{S}(i-1)) \mid R_i]$, instead of $\mathbb{E}[\log P(\mathcal{S}(i) \mid \mathcal{S}(i-1))]$, which can easily be seen to be finite. Similarly, the last term $\log P(\mathcal{T} \mid \mathcal{S}(0), \ldots, \mathcal{S}(M(L_t)))$ also depends on $R_{M(L_t)}$ through $U$, and hence we need to consider $\mathbb{E}[\log P(\mathcal{T} \mid \mathcal{S}(0), \ldots, \mathcal{S}(M(L_t))) \mid U, L_t]$.

To this end, it is instrumental to understand the inter-dependence structure among $\{\Delta_i\}_{i \geq 1}$, $\{\rho_i\}_{i \geq 0}$ and $\{R_i\}_{i \geq 1}$.

LEMMA 3.4. *The following hold*:

1. *Both*

$$\{\rho_j\}_{j\geq 1} \quad \text{under } \mathbb{P} \quad \text{and}$$

$$\{(\Delta_{\rho_j+k})_{k=1,\ldots,\rho_{j+1}-\rho_j} : j \geq 1\} \quad \text{under } \mathbb{P}(\cdot \mid \rho_j : j \geq 1)$$

*are independent families.*

2. *The $\rho_{j+1} - \rho_j$ ($j \geq 1$) has the same law as $\rho_1$, which is given by*

$$\mathbb{P}(\rho_1 = k) = \frac{k-1}{k!} \quad \text{for all } k \geq 2.$$

*Moreover, conditioned on $\{\rho_1 = k\}$, $R_1 - R_0$ is Gamma distributed with parameter $(k, 7)$. That is, it has the probability density*

$$\frac{7^k}{(k+1)!} r^{k-1} e^{-7r} \mathbb{1}\{r \geq 0\}.$$

3. *Let $\{E_i\}_{i\in\mathbb{N}}$ be independent exponential random variables with rate 7. Conditioned on $\{\rho_1 = k, T_k - T_0 = s\}$,*

$$\sum_{i=\{1,\ldots,k\}\setminus\{k-1\}} \Delta_i^{-1} \stackrel{d}{=} \frac{1}{s} \sum_{i=2}^{k} \frac{\sum_{j=1}^{k} E_j}{\sum_{j=1}^{i} \frac{1}{k-j} E_j}.$$

PROOF. The first assertion follows from the fact that $(\rho_j)_{j\geq 1}$ are stopping times for the process $(T_i)_{i\geq 0}$.

To prove the second and third assertions, it is useful to realize the interarrival times in such a way that the dependence structure between $\rho_1$, $T_k - T_0 = \sum_{i=1}^{k} \Delta_i$ and $\Delta_i^{-1}$ is clear. To this end, let $(\Delta_i^{(k)})_{i=1}^{k}$ be an increasing order statistics of independent Exp(7) random variables and let $\pi$ be a uniform random variable on the permutations $\mathfrak{S}_k$ over $\{1, 2, \ldots, k\}$, which is independent of $\Delta^{(k)}$. Then we can realize the interarrival times as

$$(3.6) \qquad (\Delta_i)_{1\leq i\leq k} = (\Delta_{\pi(i)}^{(k)})_{1\leq i\leq k}.$$

Now, since $\{\rho_1 = k\}$ depends only on $\pi$, we find

$$\mathbb{P}(\rho_1 = k) = \mathbb{P}(\Delta_1 > \Delta_2 > \cdots > \Delta_{k-1} \text{ and } \Delta_{k-1} < \Delta_k) = \frac{k-1}{k!}$$

by simply counting the number of permutations satisfying the above ordering. For the same reason, $\{\rho_1 = k\}$ is independent of $\sum_{i=1}^{k} \Delta_i = \sum_{i=1}^{k} \Delta_i^{(k)}$, which is Gamma distributed with parameter $(k, 7)$. Thus the second assertion is proved.

Finally, $\sum_{i=1}^{k} \Delta_i$ is independent of $\{\Delta_j / \sum_{i=1}^{k} \Delta_i\}_{j=1}^{k}$; see [12], Chapter IX, Theorem 4.1. Therefore, conditioned on $\{\rho_1 = k, \sum_{i=1}^{k} \Delta_i = s\}$, we have

$$(3.7) \qquad \sum_{i\in\{1,\ldots,k\}\setminus\{k-1\}} \Delta_i^{-1} \stackrel{d}{=} \frac{1}{s} \sum_{i=2}^{k} \left(\frac{\tilde{\Delta}_i^{(k)}}{\sum_{i=1}^{k} \tilde{\Delta}_i^{(k)}}\right)^{-1},$$

where $\tilde{\Delta}^{(k)}$ is an independent copy of $\Delta^{(k)}$. The third assertion follows from the following distributional identity proved in [16], §1:

$$(\tilde{\Delta}_1^{(k)}, \tilde{\Delta}_2^{(k)}, \ldots, \tilde{\Delta}_k^{(k)}) \stackrel{d}{=} \left( \frac{E_1}{k-1}, \sum_{j=1}^{2} \frac{E_j}{k-j}, \ldots, \sum_{j=1}^{k} \frac{E_j}{k-j} \right). \qquad \square$$

Now we state the bounds on the conditional expectations mentioned before.

LEMMA 3.5.

(a) *There exists $C > 0$ such that almost surely,*

$$\mathbb{E}\big[\log P\big(\mathcal{S}(1) \mid \mathcal{S}(0)\big) \mid \rho_1, R_1\big] \geq -C\left(R_1 + \frac{\rho_1^3}{R_1}\right), \tag{3.8}$$

*and*

$$\mathbb{E}\big[\log P\big(\mathcal{T} \mid \mathcal{S}(0), \ldots, \mathcal{S}(M(L_t))\big) \mid U, \sigma, L_t\big] \mathbb{1}\{U > 0\}$$
$$\geq -C\left(U + \frac{\sigma^3}{U} + (t - L_t) + (t - L_t)^{-1}\right). \tag{3.9}$$

(b) *There exists $C > 0$ such that almost surely,*

$$\mathbb{E}\big[\log P\big(\mathcal{S}(1) \mid \mathcal{S}(0)\big) \mid R_1\big] \geq -C(R_1 + R_1^{-1}), \tag{3.10}$$

*and*

$$\mathbb{E}\big[\log P\big(\mathcal{T} \mid \mathcal{S}(0), \ldots, \mathcal{S}(M(L_t))\big) \mid U, L_t\big] \mathbb{1}\{U > 0\}$$
$$\geq -C(U + U^{-1} + (t - L_t) + (t - L_t)^{-1}). \tag{3.11}$$

PROOF. *Part (a)*: By (3.4) and Lemma 3.4, we get for $\omega \in \{\rho_1 = n+2, T_{n+2} - T_0 = s\}$,

$$\log P(\mathcal{S}_1 \mid \mathcal{S}_0) \geq -C\left(s + \sum_{\substack{i=1,\ldots,n+2 \\ i \neq n+1}} \frac{1}{\Delta_i}\right) \tag{3.12}$$
$$\stackrel{d}{=} -C\left(s + \frac{1}{s} \sum_{i=2}^{n+2} \frac{\sum_{j=1}^{n+2} E_j}{\sum_{j=1}^{i} \frac{1}{n+2-j} E_j}\right).$$

Thus it suffices to show that the expectation over $\{E_1, \ldots, E_{n+2}\}$ in the last line is bounded by $(n+1)^3$. To this end, we first bound the expectation of the sum as

follows:

$$\mathbb{E}\left[\sum_{i=2}^{n+2} \frac{\sum_{j=1}^{n+2} E_j}{\sum_{j=1}^{i} \frac{1}{n+2-j} E_j}\right]$$

(3.13)
$$\leq \sum_{i=2}^{n+2} (n+2-i) \mathbb{E}\left[\frac{\sum_{j=1}^{n+2} E_j}{\sum_{j=1}^{i} E_j}\right]$$

$$= \sum_{i=2}^{n+2} (n+2-i)\left(1 + \mathbb{E}\left[\sum_{j=i+1}^{n+2} E_j\right] \mathbb{E}\left[\left(\sum_{j=1}^{i} E_j\right)^{-1}\right]\right).$$

Now this is the point where we use the extra integrability brought by omitting $i = 1$, which corresponds to the largest value of $\{\Delta_i^{-1}\}_{i=1}^n$. Indeed, since $\sum_{j=1}^{i} E_j$ is Gamma distributed with parameters $(i, 1)$, for $i \geq 2$, we can compute

$$\mathbb{E}\left[\sum_{j=i+1}^{n+2} E_j\right] = n+2-i \quad \text{and} \quad \mathbb{E}\left[\left(\sum_{j=1}^{i} E_j\right)^{-1}\right] = \frac{1}{i-1}.$$

Substituting these into (3.13), we arrive at

$$\mathbb{E}\left[\sum_{i=2}^{n+2} \frac{\sum_{j=1}^{n+2} E_j}{\sum_{j=1}^{i} \frac{1}{n+2-j} E_j}\right] \leq n \sum_{i=2}^{n+2} \frac{n+1}{i-1} \leq (n+1)^3.$$

The proof of (3.9) is essentially the same. We assume $U > 0$ and $\sigma = n$. Then recall that by (3), we have

$$\log P(\mathcal{T} \mid \mathcal{S}(0), \ldots, \mathcal{S}(M(L_t)))$$

$$\geq -C\left(\sum_{i=1}^{n-1} \Delta_i^{-1} + U + (t - L_t) + (t - L_t)^{-1}\right).$$

Since the interarrival times of disasters in $[R_{M(L_t)}, L_t]$ are decreasing, the largest member of $\{\Delta_i^{-1}\}_{i=1}^n$ is omitted in the sum on the right-hand side. This is the same situation as in Lemma 3.4(3), and thus conditioned on $U$, we have

$$\sum_{i=1,\ldots,n-1} \Delta_i^{-1} \stackrel{\mathrm{d}}{=} U^{-1} \sum_{i=2}^{n} \frac{\sum_{j=1}^{n} E_j}{\sum_{j=1}^{i} \frac{1}{n-j} E_j}.$$

Then the same computation as in the previous case yields the desired bound.

*Part (b)*: In order to take an expectation over $\rho_1$ conditioned on $R_1$, we estimate the conditional probability

$$\mathbb{P}(\rho_1 = n+2 \mid R_1 = r) = \mathbb{P}(\rho_1 = n+2 \mid T_{\rho_1} - T_0 = r)$$

$$= \frac{\mathbb{P}(\rho_1 = n+2, T_{n+2} - T_0 = r)}{\mathbb{P}(T_{\rho_1} - T_0 = r)},$$

where here and in what follows, conditions like $T_{n+2} - T_0 = r$ should be understood in the sense of probability density. Since $\{\rho_1 = n+2\}$ and $T_{n+2} - T_0$ are independent, by using Lemma 3.4, we can bound the numerator from above by

$$(3.14) \quad \mathbb{P}(\rho_1 = n+2, T_{n+2} - T_0 = r) \leq \frac{(n+1)}{(n+2)!} \frac{(7r)^{n+1}}{(n+1)!} e^{-7r}.$$

On the other hand, the denominator is bounded from below by

$$\begin{aligned}
\mathbb{P}(\rho_1 = 2, T_2 - T_0 = r) \\
= \mathbb{P}(T_1 - T_0 < T_2 - T_1, T_2 - T_0 = r) \\
(3.15) \quad = \frac{1}{2} \mathbb{P}(T_2 - T_0 = r) \\
= \frac{3}{2} r e^{-7r}.
\end{aligned}$$

Combining (3.14) and (3.15), we find the bound

$$\mathbb{P}(\rho_1 = n+2 \mid R_1 = r) \leq \frac{(n+1)}{(n+2)!} \frac{(7r)^{n+1}}{(n+1)!} \frac{2}{3r}$$

$$\leq 4 \frac{(7r)^n}{(n!)^2}.$$

In particular, we get that if $R_1 \leq \frac{1}{7}$ then

$$\mathbb{P}(\rho_1 = n+2 \mid R_1) \leq \frac{4}{(n!)^2}$$

and consequently,

$$(3.16) \quad \begin{aligned}
\mathbb{E}[\rho_1^3 \mid R_1] &= \sum_{n=0}^{\infty} (n+2)^3 \mathbb{P}(\rho_1 = n+2 \mid R_1) \\
&\leq \sum_{n=0}^{\infty} 4 \frac{(n+2)^3}{(n!)^2} < \infty.
\end{aligned}$$

If $R_1 > \frac{1}{7}$, then we use $n! \geq (\frac{n}{2})^{\frac{n}{2}}$ to see that for all $n > \sqrt{28 R_1}$, we have

$$\mathbb{P}(\rho_1 = n+2 \mid R_1) \leq \frac{4}{n^3 2^n}$$

and consequently,

$$(3.17) \quad \mathbb{E}[\rho_1^3 \mid R_1] \leq 28^2 R_1^2 + 4 \sum_{n > \sqrt{28 R_1}} \frac{(n+2)^3}{n^3 2^n}.$$

Since the sum on the right-hand side converges, we can combine the two estimates (3.16) and (3.17) to find $C > 0$ such that for all $R_1 > 0$,

$$\mathbb{E}[\rho^3 \mid R_1] \leq C(1 + R_1^2).$$

Plugging this in (3.8), we get (3.10).

Finally (3.11) follows in a similar way. We consider the probability of $\{\sigma = n\}$ conditioned on $\{U = u, L_t = l\}$, which can be written as

$$\mathbb{P}(\sigma = n \mid U = u, L_t = l)$$
$$= \frac{\mathbb{P}(\sum_{i=M(l)+1}^{M(l)+n+1} \Delta_i = u, \Delta_{M(l)+1} > \cdots > \Delta_{M(l)+n+1})}{\mathbb{P}(\sum_{i=M(l)+1}^{M(l)+\sigma+1} \Delta_i = u, \Delta_{M(l)+1} > \cdots > \Delta_{M(l)+\sigma+1})}.$$

The two events in the numerator are independent, and hence the numerator is bounded (in the sense of density) from above by

$$\frac{1}{(n+1)!} \frac{1}{n!} (7u)^n e^{-7u}. \tag{3.18}$$

On the other hand, the denominator is bounded from below by considering the special case $\sigma = 0$:

$$\mathbb{P}\left(\sum_{i=M(l)+1}^{M(l)+\sigma+1} \Delta_i = u, \Delta_{M(l)+1} > \cdots > \Delta_{M(l)+\sigma+1}\right)$$
$$\geq \mathbb{P}(\Delta_{M(l)+1} = u) \tag{3.19}$$
$$= 7e^{-7u}.$$

From (3.18) and (3.19), we find that

$$\mathbb{P}(\sigma = n \mid U = u, L_t = l) \leq \frac{1}{n+1} \frac{(7u)^n}{(n!)^2}.$$

The rest of the argument is the same as for (3.10). □

We are now ready to prove Lemma 3.1.

PROOF OF LEMMA 3.1. *Part (i)*: Note that on the event $\{M(t) = m\}$, we have

$$\log P(\mathcal{S}_t) = \log P(\mathcal{S}(0)) + \sum_{i=1}^{m} \log P(\mathcal{S}(i) \mid \mathcal{S}(i-1))$$
$$+ \log P(\mathcal{T} \mid \mathcal{S}(0), \ldots, \mathcal{S}(m)).$$

By using the bounds (3.10) and (3.11) and denoting $R_i - R_{i-1}$ by $\Delta R_i$, we get on $\{F_t < t\}$

$$\mathbb{E}[\log P(\mathcal{S}_t) \mid F_t, L_t]$$

(3.20)
$$\geq -C\left(F_t + F_t^{-1} + \mathbb{E}\left[\sum_{i=1}^{M(t)} \Delta R_i + U\right]\right.$$
$$+ \mathbb{E}\left[\sum_{i=1}^{M(t)} (\Delta R_i)^{-1} + U^{-1}\right] + (t - L_t) + (t - L_t)^{-1}\right).$$

Since we have $F_t + \sum_{i=1}^{M(t)} \Delta R_i + U + (t - L_t) = t$ by definition, it remains to show that the third expectation in (3.20) is bounded by $Ct$. We use that $A_i' \preceq_{\text{st}} A_i \preceq_{\text{st}} \Delta R_i$, where $A_i$ is Gamma distributed with parameter $(2, 7)$ and $A_i'$ is exponentially distributed with parameter 7, respectively. Since

$$(r_1, \ldots, r_i) \mapsto \frac{1}{r_1} \mathbb{P}(r_1 + \cdots + r_i \leq t)$$

is decreasing, the above stochastic domination implies

$$\mathbb{E}\left[\sum_{i=1}^{M(t)} (\Delta R_i)^{-1}\right] = \sum_{i=1}^{\infty} \mathbb{E}[(\Delta R_1)^{-1} \mathbb{1}\{\Delta R_1 + \cdots + \Delta R_i \leq t\}]$$

$$\leq \sum_{i=1}^{\infty} \mathbb{E}[A_1^{-1} \mathbb{1}\{A_1 + A_2' + \cdots + A_i' \leq t\}].$$

By using the form of the probability density of $A_1$, we find

$$\mathbb{E}[A_1^{-1} \mathbb{1}\{A_1 + A_2' + \cdots + A_i' \leq t\}]$$
$$= \int_0^\infty a^{-1} \mathbb{P}(a + A_2' + \cdots + A_i' \leq t) 49 a e^{-7a} \, da$$
$$= 7 \mathbb{P}(A_1' + \cdots + A_i' \leq t)$$

and hence

$$\mathbb{E}\left[\sum_{i=1}^{M(t)} (\Delta R_i)^{-1}\right] = 7 \sum_{i=1}^{\infty} \mathbb{P}(A_1' + \cdots + A_i' \leq t).$$

The sum on the right-hand side is nothing but the expectation of a Poisson process with intensity 7 on $[0, t]$, which is equal to $7t$.

*Part (ii)*: We follow the same strategy as in (i) but we skip (S6) in our strategy. Then we obtain the bound

$$\mathbb{E}[\log P(\tau_\infty^1(\omega) \geq t, B(s) \in J^{(6)} \text{ for all } 0 \leq s \leq t) \mid F_t]$$
$$\geq -C(t + \mathbb{1}\{F_t < t\} F_t^{-1}).$$

Since $F_t(\omega_{[0,1]^c}) \geq 1$, we are done. □

**4. Higher moments with general endpoints distribution.** In this section, we use Lemma 3.1 to get bounds on higher moments for the survival probability with more general initial and terminal distribution for the Brownian bridge. We first introduce some more notation. Given $0 \le r < s$ and $\nu \in \mathcal{M}(\mathbb{R}^2)$, we denote by $P^{\nu,r,s}$ the law of the the Brownian bridge in the interval $[r, s]$ with initial and terminal points chosen according to $\nu$. More precisely, let us recall that $P^{r,x;s,y}$ denotes the Brownian bridge measure between $(r, x)$ and $(s, y)$ and define

$$(4.1) \qquad P^{\nu,r,s}(\cdot) := \int_{\mathbb{R}^2} P^{r,x;s,y}(\cdot)\nu\bigl(\mathrm{d}(x, y)\bigr).$$

As we mentioned in Section 2, we will derive our moment bound by considering the survival probabilities in many disjoint tubes. For $x \in \mathbb{R}$ and $i \ge 1$, we define

$$J_x^{(5)}(i) := x + 7i + \left[-\frac{5}{2}, \frac{5}{2}\right] \subseteq \mathbb{R},$$

$$J_x^{(6)}(i) := x + 7i + [-3, 3] \subseteq \mathbb{R},$$

and for a given probability measure $\nu \in \mathcal{M}(\mathbb{R}^2)$ and $p \ge 1$,

$$(4.2) \qquad M^p(\nu) := \sup_{x,y \in \mathbb{R}} \min_{i=0,\ldots,p} \nu\bigl(J_x^{(5)}(i) \times J_y^{(5)}(i)\bigr).$$

This is a measure of (local) dispersion of $\nu$. If $M^p(\nu)$ is large, then under $P^{\nu,r,s}$, there is a good chance to find the initial and terminal points of the Brownian motion in $J_x^{(5)}(i) \times J_y^{(5)}(i)$, for each $i = 0, 1, \ldots, p$. Note that from our choice of $J^{(7)}$, the tubes connecting $J_x^{(5)}(i)$ and $J_y^{(5)}(i)$ are independent for different $i$. Since we can apply Lemma 3.1 to get a lower bound on the survival probability for each tube by itself, we should be able to get a better bound on the survival probability in the time interval $[r, s]$.

The following lemma is the goal of this section, which formalizes the above intuition.

LEMMA 4.1. *For every $p \ge 1$, there exists $C > 0$ such that for any $0 \le r < s$, $t \in [0, s-r]$ and $\nu \in \mathcal{M}(\mathbb{R}^2)$,*

$$(4.3) \quad \mathbb{E}\bigl[\bigl|\log P^{\nu,r,s}\bigl(\tau_\infty(\omega_{[r,s]}) \ge r + t\bigr)\bigr|^p\bigr] \le C\bigl(1 + t^p\bigr) + \bigl|\log M^p(\nu)\bigr|^p.$$

*If in addition $\nu$ is supported on $[-A, A]^2 \subseteq \mathbb{R}^2$ for some $A \ge \frac{7(3+p)}{2}$, then*

$$(4.4) \quad \begin{aligned} &\mathbb{E}\bigl[\bigl|\log P^{\nu,r,s}\bigl(\tau_\infty(\omega_{[r,s]}) \ge r + t, |B(u)| \le A \text{ for all } u \in [r, r+t]\bigr)\bigr|^p\bigr] \\ &\le C\bigl(1 + t^p\bigr) + \bigl|\log M^{p+2}(\nu)\bigr|^p. \end{aligned}$$

PROOF. We assume that the supremum in (4.2) is attained at $x, y \in \mathbb{R}$, and set

$$\nu_i\bigl(\mathrm{d}(x, y)\bigr) := \frac{\nu|_{J_x^{(5)}(i) \times J_y^{(5)}(i)}(\mathrm{d}(x, y))}{\nu(J_x^{(5)}(i) \times J_y^{(5)}(i))}.$$

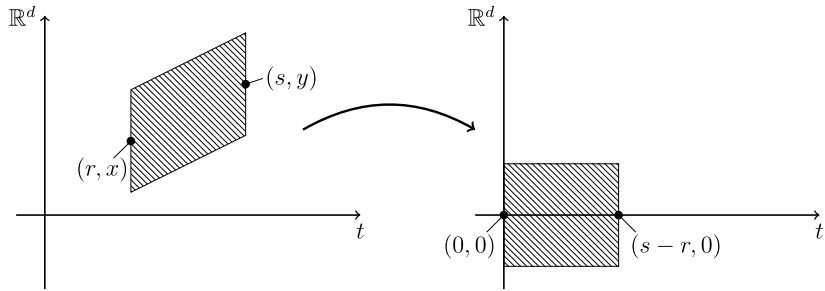

FIG. 3. *The law of $\omega$ is invariant under the affine transformation that maps $(r, x)$ to $(0, 0)$ and $(s, y)$ to $(r - s, 0)$. Note that the shifted tube connecting $\{r\} \times (x + [-3, 3])$ and $\{s\} \times (y + [-3, 3])$ is mapped onto $[0, s - r] \times J^{(6)}$, the tube considered in Lemma* 3.1.

Then we have

$$P^{\nu,r,s}(\tau_\infty(\omega_{[r,s]}) \geq r + t)$$
$$\geq \max_{i=0,\ldots,p} \int_{(u,v) \in J_x^{(5)}(i) \times J_y^{(5)}(i)} P^{r,u;s,v}(\tau_\infty(\omega_{[r,s]}) \geq r + t) \nu(\mathrm{d}(u, v))$$
$$= \max_{i=0,\ldots,p} \nu(J_x^{(5)}(i) \times J_y^{(5)}(i)) P^{\nu_i, r, s}(\tau_\infty(\omega_{[r,s]}) \geq r + t)$$
$$\geq \min_{i=0,\ldots,p} \nu(J_x^{(5)}(i) \times J_y^{(5)}(i)) \max_{i=0,\ldots,p} P^{\nu_i, r, s}(\tau_\infty(\omega_{[r,s]}) \geq r + t).$$

In order to apply Lemma 3.1 to the probability in the last line, we perform a time-space affine transformation that maps $(r, x)$ to $(0, 0)$ and $(s, y)$ to $(s - r, 0)$ (see Figure 3), and write $\bar{\omega}$ for the image of $\omega$ and $\bar{\nu}_i \in \mathcal{M}(J_0^{(5)}(i)^2)$ for the image measure of $\nu_i$, respectively. Under this transformation, $\bar{\omega}$ has the same law as $\omega$ while $P^{\nu_i,r,s}$ is transformed to $P^{\bar{\nu}_i} = P^{\bar{\nu}_i, 0, s-r}$. Therefore, we have

(4.5)
$$P^{\nu_i, r, s}(\tau_\infty(\omega_{[r,s]}) \geq r + t)$$
$$\geq P^{\bar{\nu}_i}(\tau_\infty(\bar{\omega}) \geq t, B(u) \in J_0^{(6)}(i) \text{ for all } 0 \leq u \leq t)$$

and for different $i$'s, the probabilities on the right-hand side depend on $\bar{\omega}$ in disjoint sets and hence are independent under $\mathbb{P}$. Let us introduce

(4.6) $\quad X_i := \left| \log P^{\bar{\nu}_i}(\tau_\infty(\bar{\omega}) \geq t, B(u) \in J_0^{(6)}(i) \text{ for all } 0 \leq u \leq t) \right|,$

so that we can write

$$\left|\log P^{\nu,r,s}(\tau_\infty(\omega_{[r,s]}) \geq r + t)\right|^p \leq 2^{p-1} \left( |\log M^p(\nu)|^p + \left(\min_{i=0,\ldots,p} X_i\right)^p \right).$$

It remains to bound the $p$-th moment of $\min_{i=0,\ldots,p} X_i$. Recall that $\bar{\omega}$ has the same law as $\omega$ and that $X_i = X_0 \circ \theta_{0,7i}$, where $\theta_{0,7i}$ is the time-space shift operator. To simplify the notation, we write $F_{i,t}$ and $L_{i,t}$ for $F_t \circ \theta_{0,7i}(\bar{\omega})$ and $L_t \circ \theta_{0,7i}(\bar{\omega})$, respectively. By Lemma 3.1, we have the following upper bound for $i = 0, 1, \ldots, p$:

$$\mathbb{E}[X_i \mid F_{i,t}, L_{i,t}] \leq c_1 \left( t + \mathbb{1}\{F_{i,t} < t\} \left( \frac{1}{F_{i,t}} + \frac{1}{t - L_{i,t}} \right) \right),$$

where $c_1 > 0$ is a constant. In this proof, we keep the constants indexed and clarify their dependence on parameters. Using Jensen's inequality, the above bound and that the marginal laws of $F_{i,t}$ and $L_{i,t}$ are the exponential law with rate 7 truncated at $t$, we obtain that for any $\epsilon > 0$, there exists $c_2(\epsilon) > 0$ such that for $i = 0, 1, \ldots, p$,

$$\mathbb{E}\left[ (X_i - c_1 t)_+^{1-\epsilon} \right] \leq c_1 \mathbb{E}\left[ \mathbb{1}\{F_{i,t} < t\} \left( \frac{1}{F_{i,t}} + \frac{1}{t - L_{i,t}} \right)^{1-\epsilon} \right]$$

$$\leq c_2(\epsilon).$$

This bound and the Markov inequality yield

$$\mathbb{P}\left( (X_i - c_1 t)_+ \geq u \right) \leq c_2(\epsilon) u^{\epsilon - 1}$$

for all $i = 0, 1, \ldots, p$ and $u > 0$. As a consequence, if we choose $\epsilon$ sufficiently small, we have

$$\mathbb{P}\left( \min_{0 \leq i \leq p} (X_i - c_1 t)_+ \geq u \right) = \prod_{i=0}^{p} \mathbb{P}\left( (X_i - c_1 t)_+ \geq u \right)$$

$$\leq c_2(\epsilon)^{p+1} u^{-p-1/2}$$

for all $u > 0$, where we used that $X_0, X_1, \ldots, X_p$ are independent. From this tail bound, we can deduce that

$$\mathbb{E}\left[ \left( \min_{i=0,\ldots,p} X_i \right)^p \right] \leq \mathbb{E}\left[ \left( \min_{i=0,\ldots,p} (X_i - c_1 t)_+ + c_1 t \right)^p \right]$$

$$\leq c_3(p) \left( t^p + \mathbb{E}\left[ \left( \min_{i=0,\ldots,p} (X_i - c_1 t)_+ \right)^p \right] \right)$$

$$\leq c_3(p) \left( t^p + \int_0^\infty p u^{p-1} c_2(\epsilon)^{p+1} u^{-p-1/2} du \right)$$

$$\leq c_4(p, \epsilon)(t^p + 1).$$

This completes the proof of the first assertion.

The second assertion is essentially proved in the above argument once we account for some issues with the boundary. Note that the bound (4.4) is trivial unless $M^{p+2}(\nu) > 0$, and in that case we again write $(x, y)$ for the values where the supremum in (4.2) is attained. Since $A$ is large enough, we observe that among

$$\{J_x^{(5)}(i) \times J_y^{(5)}(i) : i = 0, \ldots, p + 2\}$$

there are at least $p + 1$ indices $i_0, \ldots, i_p$ such that
$$J_x^{(6)}(i_j) \times J_y^{(6)}(i_j) \subseteq [-A, A]^2 \quad \text{for all } j = 0, \ldots, p.$$

For such an index $i_j$, we note that the event considered in (4.6) ensures that the Brownian motion does not leave $[-A, A]$ in $[r, r+t]$. We then obtain (4.4) by the same argument as for the first assertion where $\min_{i=0,\ldots,p} X_i$ has to be replaced by $\min_{j=0,\ldots,p} X_{i_j}$. □

**5. Midpoint distribution of polymer.** In order to prove Proposition 2.7, we will apply Lemma 4.1 to the midpoints distribution under the following polymer measures

(5.1)
$$\nu_{\omega,\beta}^{r,s,t}(\mathrm{d}(x, y))$$
$$:= P\big((B(r), B(s)) \in \mathrm{d}(x, y) \mid \tau_\beta^1(\omega_{[r,s]^c}) \geq t, \mathcal{A}_t\big) \in \mathcal{M}(\mathbb{R}^2).$$

Thus we need to estimate the dispersion $M^p$ of this measure, which is the goal of this section:

LEMMA 5.1. *Let $p \geq 0$ and $q \geq 1$. There exists $C > 0$ such that for all $\beta \in [0, \infty]$ and all $1 \leq r^- \leq r^+ \leq t$ such that either $r^+ \leq t - 1$ or $r^+ = t$,*

(5.2)
$$\mathbb{E}\big[|\log M^p(\nu_{\omega,\beta}^{r^-,r^+,t})|^q\big] \leq C(1 + \log^+ t)^C.$$

PROOF. Let us recall the notation
$$J_x^{(1)} = x + \left[\frac{1}{2}, \frac{1}{2}\right),$$
$$J_x^{(5)}(i) = x + 7i + \left[-\frac{5}{2}, \frac{5}{2}\right],$$
$$J_x^{(6)}(i) = x + 7i + [-3, 3],$$
$$M^p(\nu) = \sup_{x, y \in \mathbb{R}} \min_{i=0,\ldots,p} \nu\big(J_x^{(5)}(i) \times J_y^{(5)}(i)\big).$$

Observe first that thanks to the truncation $\mathcal{A}_t$, for every $0 \leq r < s \leq t$ and every $\omega$, there exist $x, y \in \mathbb{R}$ such that

(5.3) $$\nu_{\omega,\beta}^{r,s,t}\big(J_x^{(1)} \times J_y^{(1)}\big) \geq ct^{-4}.$$

The bound (5.2) for $p = 0$ follows by setting $r = r^-$ and $s = r^+$.

In order to prove (5.2) for $p \geq 1$, we need to find sets of intervals $\{J_x^{(5)}(i)\}_{i=0}^p$ and $\{J_y^{(5)}(i)\}_{i=0}^p$ for which $\nu_{\omega,\beta}^{r^-,r^+,t}(J_x^{(5)}(i) \times J_y^{(5)}(i))$ are not too small for all $i \in \{0, 1, \ldots, p\}$. Our strategy is to use (5.3) for some $r < r^-$ and $s > r^+$ first and then *sprinkle* the mass on the time intervals $[r, r^-]$ and $[r^+, s]$. To this end, we have to find $r < r^- < r^+ < s$ and $x, y \in \mathbb{R}$ such that

- (5.3) is satisfied,
- there are no obstacles inside $[r, r^-]$ and $[r^+, s]$, close to $(r, x)$ or $(s, y)$.

The latter condition would ensure that the disasters do not prevent sprinkling the mass.

For now, let us assume that $r^- \geq 2$ and $r^+ \leq t - 1$. We denote $r_0^- := r^- - 1$ and $r_0^+ := r^+ + 1$ and for $i \geq 1$,

$$(5.4) \qquad r_i^- := r_0^- + \frac{6}{\pi^2} \sum_{j=1}^{i} j^{-2} \quad \text{and} \quad r_i^+ := r_0^+ - \frac{6}{\pi^2} \sum_{j=1}^{i} j^{-2}.$$

Note that $r_i^- < r^-$ and $r_i^+ > r^+$ for all $i$. From (5.3), we know that there exists $(j_i^+, j_i^-)$ such that

$$(5.5) \qquad v_{\omega,\beta}^{r_i^-, r_i^+, t}(J_{j_i^-}^{(1)} \times J_{j_i^+}^{(1)}) \geq ct^{-4}.$$

For $i \geq 0$, let $\lambda_i : [0, \infty) \to \mathbb{R}$ be the affine linear function with $\lambda_i(r_i^-) = j_i^-$ and $\lambda_i(r_i^+) = j_i^+$, and introduce the slanted space-time boxes

$$(5.6) \qquad \begin{aligned} S_i^{\pm} := \Big\{ (u, x) : u \in [r_i^{\pm}, r_{i+1}^{\pm}), \lambda_i(u) - \frac{7}{2} \\ \leq x \leq \lambda_i(u) + 7(p+1)(q+1) - \frac{7}{2} \Big\}. \end{aligned}$$

Here, we interpret the time-interval $[r_i^+, r_{i+1}^+)$ as $(r_{i+1}^+, r_i^+]$ by a slight abuse of notation. The same convention applies in the rest of this proof. Let us define the event

$$C_i := \{\omega(S_i^+ \cup S_i^-) = 0\}.$$

Observe that since the boxes $S_i^{\pm}$ are disjoint and have decreasing volume, the events are independent and $\mathbb{P}(C_i) \geq \mathbb{P}(C_0) > 0$ for all $i \geq 0$. Therefore,

$$G := \inf\{i \geq 0 : C_i \text{ holds}\}$$

has a geometric tail:

$$(5.7) \qquad \mathbb{P}(G \geq i) \leq (1 - \mathbb{P}(C_0))^i.$$

In particular $G$ is almost surely finite, and hence $j_G^-$ and $j_G^+$ are well defined. Figure 4 provides a schematic picture of this construction.

Now for $k \in \{0, \ldots, q\}$, $l \in \{0, \ldots, p\}$ and $u \geq 0$, let

$$J^{(5)}(k, l, u) := J^{(5)}_{\lambda_G(u)+7(p+1)k}(l),$$

$$J^{(6)}(k, l, u) := J^{(6)}_{\lambda_G(u)+7(p+1)k}(l),$$

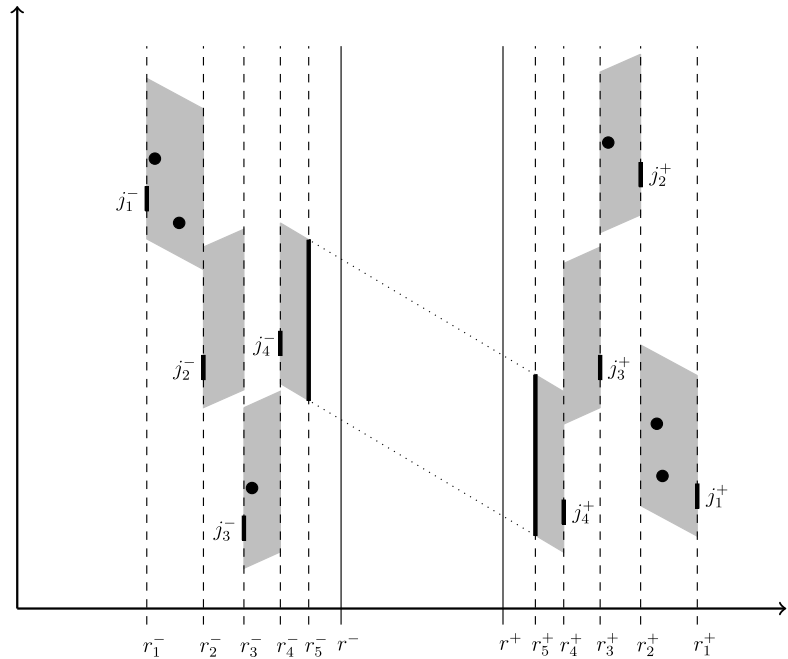

FIG. 4. *An illustration of the resampling procedure. The dots represent the disasters. The short black intervals $J^{(1)}_{j_i^-} \times J^{(1)}_{j_i^+}$ have not too small probability ($\geq ct^{-2}$) under the polymer measure with respect to $\omega_{[r_i^-, r_i^+]^c}$. The gray areas corresponds to $S_i^\pm$. In this figure, we have $G = 4$ since the $S_4^+$ and $S_4^-$ are free of disasters. The mass of the polymer measure in $J^{(1)}_{j_4^\pm}$ can be sprinkled to the long black intervals.*

and for $\pm \in \{+, -\}$, consider space-time tubes

$$J^{(6)}_\pm(k, l) := \{(u, x) : u \in [r_{G+1}^\pm, r^\pm], x \in J^{(6)}_\pm(k, l, u)\},$$
$$J^{(6)}_\pm(k) := J^{(6)}_\pm(k, 0) \cup \cdots \cup J^{(6)}_\pm(k, p).$$

We define the events

$$\mathcal{A}_1^\pm(k, l) := \{(u, B(u)) \in S_G \text{ for all } u \in [r_G^\pm, r_{G+1}^\pm],$$
$$B(r_{G+1}^\pm) \in J^{(5)}_\pm(k, l, r_{G+1}^\pm)\},$$
$$\mathcal{A}_2^\pm(k, l) := \{(u, B(u)) \in J^{(6)}_\pm(k, l) \setminus \mathcal{D} \text{ for all } u \in [r_{G+1}^\pm, r^\pm],$$
$$B(r^\pm) \in J^{(5)}_\pm(k, l, r^\pm)\},$$

where $\mathcal{D}$ is the set of disasters defined in (1.6). In words, $\mathcal{A}_1^-(k, l)$ is the event that the Brownian motion moves from $J^{(1)}_{j_G^-}$ to the left end of the tube $J^{(6)}_-(k, l)$ in $[r_G^-, r_{G+1}^-]$, without leaving $S_G^-$. This guarantees survival since by the definition of $G$, there are no disasters in $S_G^-$. On the other hand, $\mathcal{A}_2^-(k, l)$ is the event that the

Brownian motion survives inside tube $J_-^{(6)}(k,l)$ in $[r_{G+1}^-, r^-]$. We set

$$\mathcal{A}(k,l) := \mathcal{A}_1^-(k,l) \cap \mathcal{A}_2^-(k,l) \cap \mathcal{A}_2^+(k,l) \cap \mathcal{A}_1^+(k,l).$$

By definition, we know that $M^p(v_{\omega,\beta}^{r^-,r^+,t})$ is bounded from below by the $\max_{k \in \{0,1,\ldots,q\}} \min_{l \in \{0,1,\ldots,p\}}$ of the following probability:

$$v_{\omega,\beta}^{r^-,r^+,t}(J^{(5)}(k,l,r^-) \times J^{(5)}(k,l,r^+))$$
$$\stackrel{\text{def}}{=} P((B(r^-), B(r^+)) \in J^{(5)}(k,l,r^-) \times J^{(5)}(k,l,r^+) \mid \tau_\beta^1(\omega_{[r^-,r^+]^c}) \geq t, \mathcal{A}_t)$$
$$\geq P((B(r_G^-), B(r_G^+)) \in J_{j_G^-}^{(1)} \times J_{j_G^+}^{(1)}, \mathcal{A}(k,l) \mid \tau_\beta^1(\omega_{[r_G^-,r_G^+]^c}) \geq t, \mathcal{A}_t),$$

where in the last line, we have used that

$$\mathcal{A}(k,l) \cap \{\tau_\beta^1(\omega_{[r_G^-,r_G^+]^c}) \geq t\} \subset \{\tau_\beta^1(\omega_{[r^-,r^+]^c}) \geq t\}.$$

Let us introduce the distribution

$$\alpha(x_1, y_1, dx_2, dy_2) := P^{r_{G+1}^-, x_1; r_{G+1}^+, y_2}((B(r^-), B(r^+)) \in d(x_2, y_2))$$

and denote

$$p(x_1, y_1, dx_2, dy_2)$$
$$:= P^{r_G^-, x_1; r_G^+, y_1}(B(r_{G+1}^-) \in dx_2, B(r_{G+1}^+) \in dy_2, (u, B(u)) \in S_G^- \cup S_G^+$$
$$\text{for all } u \in [r_G^-, r_{G+1}^-] \cup [r_{G+1}^+, r_G^+]).$$

Note that $r_{G+1}^- - r_G^- = r_G^+ - r_{G+1}^+ = 6/(\pi^2(G+1)^2)$ and, therefore,

$$\inf_{(x_1,y_1) \in J_{j_G^-}^{(1)} \times J_{j_G^+}^{(1)}} \left\{ \int_{J_-^{(5)}(k,l,r_{G+1}^+) \times J_+^{(5)}(k,l,r_{G+1}^-)} p(x_1, y_1, dx_2, dy_2) \right\} \geq e^{-cG^2}.$$

Note also that since $S_G^-$ and $S_G^+$ are slanted parallel to the line connecting $(r_G^-, j_G^-)$ and $(r_G^+, j_G^+)$, we can apply an affine transformation and use invariance of the Brownian bridge to see that this estimate does not depend on the distance between

$j_G^-$ and $j_G^+$. Using the above notation and estimate, we get

$$P((B(r_G^-), B(r_G^+)) \in J_{j_G^-}^{(1)} \times J_{j_G^+}^{(1)}, \mathcal{A}(k,l) \mid \tau_\beta^1(\omega_{[r_G^-, r_G^+]^c}) \geq t, \mathcal{A}_t)$$

$$= \int_{J_{j_G^-}^{(1)} \times J_{j_G^+}^{(1)}} v_{\omega,\beta}^{r_G^-, r_G^+, t}(d(x_1, y_1)) \int_{J^{(5)}(k,l,r_{G+1}^-) \times J^{(5)}(k,l,r_{G+1}^+)} p(x_1, y_1, dx_2, dy_2)$$

$$\times \int_{J^{(5)}(k,l,r^-) \times J^{(5)}(k,l,r^+)} \alpha(x_2, y_2, dx_3, dy_3) P^{r_{G+1}^-, x_2; r^-, x_3}(\mathcal{A}_2^-(k,l))$$

$$\times P^{r^+, y_3; r_{G+1}^+, y_2}(\mathcal{A}_2^+(k,l))$$

$$\geq v_{\omega,\beta}^{r_G^-, r_G^+, t}(J_{j_G^-}^{(1)} \times J_{j_G^+}^{(1)}) e^{-cG^2}$$

$$\times \inf_{x_2, x_3, y_2, y_3} P^{r_{G+1}^-, x_2; r^-, x_3}(\mathcal{A}_2^-(k,l)) P^{r^+, y_3; r_{G+1}^+, y_2}(\mathcal{A}_2^+(k,l)),$$

where the infimum is over $J^{(5)}(k,l,r_{G+1}^-) \times J^{(5)}(k,l,r^-) \times J^{(5)}(k,l,r^+) \times J^{(5)}(k,l,r_{G+1}^+)$. Recalling (5.5) and noting that $G$ has all moments by (5.7), we only need to prove that $\min_{k \in \{0,1,\ldots,q\}} \max_{l \in \{0,1,\ldots,p\}}$ of

$$Z_{k,l} := \left| \log \inf_{(x_2, x_3, y_2, y_3)} P^{r_{G+1}^-, x_2; r^-, x_3}(\mathcal{A}_2^-(k,l)) P^{r^+, y_3; r_{G+1}^+, y_2}(\mathcal{A}_2^+(k,l)) \right|$$

has all moments. Now letting $F_k^\pm$ and $L_k^\pm$ denote the first and last disasters in $J_\pm^{(6)}(k)$, respectively, we get from (3.1) that there exists $C > 0$ such that

$$\mathbb{E}\left[\max_{l \in \{0,1,\ldots,p\}} Z_{k,l} \mid \omega_{[r_{G+1}^-, r_{G+1}^+]^c}, F_k^+, F_k^-, L_k^+, L_k^-\right]$$

$$\leq C\left(1 + \mathbb{1}\{F_k^- < r^- - r_{G+1}^-\}\left(\frac{1}{F_k^-} + \frac{1}{r^- - L_k^-}\right)\right.$$

$$\left. + \mathbb{1}\{F_k^+ < r_{G+1}^+ - r^+\}\left(\frac{1}{F_k^+} + \frac{1}{r_{G+1}^+ - L_k^+}\right)\right).$$

Then we can argue exactly in the same way as in the proof of Lemma 4.1 to obtain

$$\mathbb{E}\left[\left(\min_{k \in \{0,1,\ldots,q\}} \max_{l \in \{0,1,\ldots,p\}} Z_{k,l}\right)^q\right] \leq C.$$

This completes the proof for $r^- \geq 2$ and $r^- \leq t - 1$.

In the case $r^- < 2$, we use the interval $[0, 1]$, which is free of disasters, in place of $[r_G^-, r_{G+1}^-]$, and set $j_i^- = 0$ for all $i \geq 0$. More precisely, define $r_i^+$ as above and let $j_i^+$ be such that

$$v_{\omega,\beta}^{1, r_i^+, t}(\mathbb{R} \times J_{j_i^+}^{(1)}) \geq Ct^{-2}.$$

Let $\lambda_i$ be the linear function with $\lambda_i(0) = 0$ and $\lambda_i(r_i^+) = j_i^+$ and define $S_i^+$ as in (5.6). Using an affine transformation similar to before, we see that there exists $C > 0$ (independent of $\omega$, $i$ or $j_i^+$) such that for all $y \in J_{j_i^+}^{(1)}$ and all $k = 0, \ldots, (p+1)(q+1) - 1$,

$$
\begin{aligned}
P^{0,0;r_i^+,y}&\big(B(1) \in J_{\lambda_i(1)}^{(5)}(k), B(r_{i+1}^+) \in J_{\lambda_i(r_{i+1}^+)}^{(5)}(k), (u, B(u)) \in S_i^+ \\
&\text{for all } u \in [r_{i+1}^+, r_i^+]\big) \\
(5.8) \quad = P^{0,0;r_i^+,y-j_i^+}&\big(B(1) \in J_0^{(5)}(k), B(r_{i+1}^+) \in J_0^{(5)}(k), (u, B(u)) \in \widetilde{S}_i^+ \\
&\text{for all } u \in [r_{i+1}^+, r_i^+]\big) \\
\geq C^{-1}e^{-Ci^2},
\end{aligned}
$$

where $\widetilde{S}_i := [r_{i+1}^+, r_i^+] \times [-\frac{7}{2}, 7(p+1)(q+1) - \frac{7}{2}]$. Now let $G := \inf\{i \geq 0 : S_i \cap \omega = \varnothing\}$ and note that $G$ has a geometric tail, so that in particular $j_G^+$ is well defined. By the same consideration as before, it follows that

$$
(5.9) \quad \min_{k=0,\ldots,(p+1)(q+1)-1} \nu_{\omega,\beta}^{1,r_{G+1}^+,t}\big(J_{\lambda_G(1)}^{(5)}(k) \times J_{\lambda_G(r_{G+1}^+)}^{(5)}(k)\big) \geq C^{-1}e^{-CG^2}t^{-2}.
$$

The rest of the argument is identical to before.

Finally, in the case $r^+ = t$, we simply restrict to $x = y$ in (4.2) to get

$$
M^p(\nu_{\omega,\beta}^{r^-,t,t}) \geq C \sup_{x \in \mathbb{R}} \min_{i=0,\ldots,p} P\big(B(r^-) \in J_x^{(5)}(i), B(t) \in J_x^{(5)}(i) \mid \tau_\beta^1(\omega_{[r^-,t]^c}) \geq t\big).
$$

Since we do not need to consider the survival strategy after time $t$, we can modify the previous argument by setting $r_i^+ = t$ and $j_i^+ = j_i^-$ for all $i \geq 0$ to show that

$$
\mathbb{E}\Big[\big(\log \sup_{x \in \mathbb{R}} \min_{i=0,\ldots,p} P\big(B(r^-) \in J_x^{(5)}(i), B(t) \in J_x^{(5)}(i) \mid \tau_\beta^1(\omega_{[r^-,t]^c}) \geq t\big)\big)^q\Big] \\
\leq C(1 + \log^+ t)^C. \qquad \square
$$

**6. Proof of the key propositions.** In this section, we derive all the key propositions appearing in Section 2, using Lemmas 4.1 and 5.1.

PROOF OF PROPOSITION 2.7. We write $\omega'$ for $\omega_{[r,r+s]^c}$ in this proof for simplicity of notation. We define a random probability measure $\nu(\omega')$ by

$$
(6.1) \quad \nu(\omega')(dx, dy) := P\big(B(r+s) \in dx, B(s) \in dy \mid \tau_\beta^1(\omega') \geq t, \mathcal{A}_t\big).
$$

Then we can write

$$\text{LHS of (2.6)} = \mathbb{E}\big[|\log P(\tau_\beta^1(\omega) \geq t \mid \tau_\beta^1(\omega') \geq t, \mathcal{A}_t)|^p\big]$$

$$\leq \mathbb{E}\Big[\Big|\log P^{\nu(\omega'),r,r+s}\Big(\tau_\infty^1(\omega_{[r,r+s]}) \geq t, \sup_{t' \in [r,r+s]}|B(s)| \leq t^2\Big)\Big|^p\Big].$$

Since $\nu(\omega')$ depends only on the environment outside of $[r, r+s] \times \mathbb{R}$, we may integrate $\omega_{[r,r+s]}$ conditionally on $\nu(\omega')$. Then by using (4.4) and part (i) of Lemma 5.1, we can find $C > 0$ such that the above right-hand side is bounded by

$$C(1 + s^p) + \mathbb{E}\big[|\log M^{p+2}(\nu(\omega'))|^p\big] \leq C(1 + s^p) + C(1 + \log^+ t)^C. \quad \square$$

PROOF OF PROPOSITION 2.1. We introduce a variation of the truncation event $\mathcal{A}_t$: For $s, t \geq 0$, let

$$\mathcal{A}_s^t := \{\sup\{|B(t')| : 0 \leq t' \leq s\} \leq \lceil t \rceil^2\}.$$

We define a random distribution

$$\mu(\omega)(dx) := P(B(s) \in dx \mid \tau_\beta^1(\omega) \geq s, \mathcal{A}_s^{s+t}) \in \mathcal{M}(\mathbb{R})$$

and we use $P^{\mu(\omega)}$ for the law of Brownian motion started with $B(0)$ distributed according to $\mu(\omega)$. Then we have

$$a_\beta(s+t) = \mathbb{E}[\log P(\tau_\beta^1(\omega) \geq s+t, \mathcal{A}_{s+t})]$$

$$= \mathbb{E}[\log P(\tau_\beta^1(\omega) \geq s, \mathcal{A}_s^{s+t})] + \mathbb{E}[\log P^{\mu(\omega)}(\tau_\beta(\theta_s(\omega)) \geq t, \mathcal{A}_t^{s+t})]$$

$$=: \mathbb{E}[\log P(\tau_\beta^1(\omega) \geq s, \mathcal{A}_s^{s+t})]$$

$$\quad + \mathbb{E}[\log P^{\mu(\omega)}(\tau_\beta^1(\theta_s(\omega)) \geq t, \mathcal{A}_t^{s+t})] + b(s,t)$$

$$\geq a_\beta(s) + \mathbb{E}[\log P^{\mu(\omega)}(\tau_\beta(\theta_s(\omega)) \geq t, \mathcal{A}_t)] + b(s,t),$$

where the remainder term is, by Proposition 2.7,

$$b(s,t) := \mathbb{E}[\log P(\tau_\beta^1(\omega) \geq s+t \mid \tau_\beta^1(\omega_{[s,s+1]^c}) \geq s+t, \mathcal{A}_{s+t})]$$

$$\geq -C(1 + \log^+(s+t))^C.$$

Since we have by Jensen's inequality that

$$\mathbb{E}[\log P^{\mu(\omega)}(\tau_\beta^1(\theta_s(\omega)) \geq t, \mathcal{A}_t)]$$

$$\geq \mathbb{E}\bigg[\int \log P^{\delta_x}(\tau_\beta^1(\theta_s(\omega)) \geq s, \mathcal{A}_t)\mu(\omega)(dx)\bigg]$$

$$= a_\beta(t),$$

the proof of (2.2) is completed. $\square$

PROOF OF PROPOSITION 2.2. First consider the case $t \in \mathbb{N}$. We regard $\omega$ as the sum of independent random measures: $\omega = \sum_{i \geq 0} \omega_{[i,i+1]}$ and apply a moment bound in [4] for functions of independent random variables. Let $\omega$ and $\omega'$ be two independent realizations of the environment, and for $i = 1, \ldots, t$, let

$$\omega_i := \omega_{[i,i+1]^c} + \omega'_{[i,i+1]}.$$

In other words, $\omega_i$ is obtained by re-sampling the disasters of $\omega$ in the stripe $[i-1, i) \times \mathbb{R}$. We set

$$X := \log P(\tau^1_\beta(\omega) \geq t, \mathcal{A}_t),$$

$$X_i := \log P(\tau^1_\beta(\omega_i) \geq t, \mathcal{A}_t).$$

Then Theorem 15.5 in [4] and Jensen's inequality tell us that there exists $C > 0$ depending only on $q$ such that

(6.2)
$$\mathbb{E}[|X - \mathbb{E}[X]|^{2q}]$$
$$\leq C\mathbb{E}\left[\left(\sum_{i=0}^{t-1} \mathbb{E}[((X-X_i)^+)^2 \mid \omega]\right)^q\right]$$
$$+ C\mathbb{E}\left[\left(\sum_{i=0}^{t-1} \mathbb{E}[((X-X_i)^-)^2 \mid \omega]\right)^q\right]$$
$$\leq C t^{q-1} \sum_{i=0}^{t-1} (\mathbb{E}[((X-X_i)^+)^{2q}] + \mathbb{E}[((X-X_i)^-)^{2q}]).$$

Since $(X - X_i)^+$ and $(X - X_i)^-$ have the same law, we focus on the first one. In our setting, we have

$$(X - X_i)^+ = \mathbb{1}\{X_i \leq X\}(\log P(\tau^1_\beta(\omega) \geq t, \mathcal{A}_t) - \log P(\tau^1_\beta(\omega_i) \geq t, \mathcal{A}_t))$$
$$\leq \mathbb{1}\{X_i \leq X\}(\log P(\tau^1_\beta(\omega_{[i,i+1]^c}) \geq t, \mathcal{A}_t) - \log P(\tau^1_\beta(\omega_i) \geq t, \mathcal{A}_t))$$
$$\leq |\log P(\tau^1_\beta(\omega_i) \geq t \mid \tau^1_\beta(\omega_{[i,i+1]^c}) \geq t)|.$$

Noting that the right-hand side depends only on $\omega_i$ that has the same law as $\omega$, we may apply Proposition 2.7 to find a constant $C > 0$ independent of $t$ and $\beta$ such that

$$\mathbb{E}[((X-X_i)^+)^{2q}] \leq \mathbb{E}[|\log P(\tau^1_\beta(\omega) \geq t \mid \tau^1_\beta(\omega_{[i,i+1]^c}) \geq t, \mathcal{A}_t)|^{2q}]$$
$$\leq C(1 + \log^+ t)^{Cq}.$$

Substituting this into (6.2), we obtain

$$\mathbb{E}[|X - \mathbb{E}[X]|^{2q}] \leq C t^q (1 + \log^+ t)^{Cq}$$

and the desired bound (2.3) for $t \in \mathbb{N}$ follows readily.

It remains to show that it suffices to consider the case $t \in \mathbb{N}$. By Proposition 2.7, we find $C > 0$ such that for all $t \geq 1$ and all $\beta \in [0, \infty]$,

$$\mathbb{E}[\log P(\tau_\beta(\omega) \geq t, \mathcal{A}_t)] - \mathbb{E}[\log P(\tau_\beta(\omega) \geq \lceil t \rceil, \mathcal{A}_t)] \leq C(1 + \log^+ t)^C.$$

Moreover by the same proposition with $p = 2q + 2$, we see that for $t$ sufficiently large,

$$\mathbb{P}(\log P(\tau_\beta^1(\omega) \geq t, \mathcal{A}_t)] - \log P(\tau_\beta^1(\omega) \geq \lceil t \rceil, \mathcal{A}_t) \geq t^{\frac{1}{2}})$$

(6.3) $\qquad \leq t^{-(q+1)} \mathbb{E}[(\log P(\tau_\beta^1(\omega) \geq t, \mathcal{A}_t)] - \log P(\tau_\beta^1(\omega) \geq \lceil t \rceil, \mathcal{A}_t))^{2q+2}]$

$\qquad \leq t^{-(q+\frac{1}{2})}.$

These two bounds allow us to extend (2.3) to $t \in \mathbb{R}_+$.  □

PROOF OF PROPOSITION 2.4. It follows from Proposition 2.1 and Theorem 2 in [14] that $p(\beta) := \lim_{t \to \infty} a_\beta(t)/t$ exists and

$$\frac{a_\beta(t)}{t} \leq p(\beta) + 4 \int_{2t}^\infty s^{-(2-\delta)} \, ds \leq p(\beta) + \frac{4}{2-\delta} t^{-(1-\delta)}.$$

In order to prove the other bound, we first prove that for $t \geq t_0$,

(6.4) $\qquad a_\beta(2t) \leq 2a_\beta(t) + Ct^{\frac{1}{2}+\delta}.$

To this end, we define for all $x \in [-\lceil t \rceil^2, \lceil t \rceil^2 - 1] \cap \mathbb{Z}$,

$$p_x := P(B(t) \in [x, x+1), \tau_\beta^1(\omega) \geq t, \mathcal{A}_t),$$

$$\mu_x(\omega)(dx) := P(B(t) \in dx \mid B(t) \in [x, x+1), \tau_\beta^1(\omega) \geq t, \mathcal{A}_t) \in \mathcal{M}([x, x+1)),$$

$$X_x := P^{\mu_x(\omega)}(\tau_\beta^1(\theta_t(\omega)) \geq t, \mathcal{A}_t),$$

$$Y_x := P^{\delta_x}(\tau_\beta^1(\theta_t(\omega)) \geq t, \mathcal{A}_t),$$

where as before $P^\mu$ denotes the law of Brownian motion started with $B(0)$ distributed according to $\mu$. Moreover, we consider events

$$\mathcal{B}_0 := \{Y_0 = \max\{Y_x : x \in [-t^2, t^2] \cap \mathbb{Z}\}\},$$

$$\mathcal{B}_1 := \{|\log P(\tau_\beta^1(\omega) \geq 2t, \mathcal{A}_{2t}) - \mathbb{E}[\log P(\tau_\beta^1(\omega) \geq 2t, \mathcal{A}_{2t})]| \leq (2t)^{\frac{1}{2}+\delta}\},$$

$$\mathcal{B}_2 := \{|\log P(\tau_\beta^1(\omega) \geq t, \mathcal{A}_t) - \mathbb{E}[\log P(\tau_\beta^1(\omega) \geq t, \mathcal{A}_t)]| \leq t^{\frac{1}{2}+\delta}\},$$

$$\mathcal{B}_3 := \{|Y_0 - \mathbb{E}[Y_0]| \leq t^{\frac{1}{2}+\delta}\}.$$

Since $\{Y_x : x \in \mathbb{Z}\}$ is a stationary sequence, we have

(6.5) $\qquad \mathbb{P}(\mathcal{B}_0) = (2t^2 + 1)^{-1}.$

Note that there exists $C > 0$ such that for all $\mu \in \mathcal{M}([0, 1])$ and all $x \in \mathbb{R}$,
$$P^{\delta_0}(B(1) \in dx) \vee P^{\delta_1}(B(1) \in dx) \geq CP^\mu(B(1) \in dx).$$
This implies that almost surely for all $x \in [-\lceil t \rceil^2, \lceil t \rceil^2 - 1] \cap \mathbb{Z}$,
$$Y_x \vee Y_{x+1} \geq CX_x.$$
By Proposition 2.3, there exists $C > 0$ such that for all $t$ and all $\beta \in [0, \infty]$
$$\mathbb{P}(\mathcal{B}_1 \cap \mathcal{B}_2 \cap \mathcal{B}_3) \geq 1 - Ct^{-3}. \tag{6.6}$$

Combining (6.5) and (6.6), we find that $\mathcal{B} := \mathcal{B}_0 \cap \mathcal{B}_1 \cap \mathcal{B}_2 \cap \mathcal{B}_4$ has positive probability for all sufficiently large $t$. In particular, it is nonempty and we can pick an $\omega \in \mathcal{B}$. Then since $\omega \in \mathcal{B}_0$, we have

$$P(\tau_\beta^1(\omega) \geq 2t, \mathcal{A}_{2t})$$
$$\leq \sum_{x \in [-t^2, t^2] \cap \mathbb{Z}} p_x X_x + P\left(\sup_{r \in [0,t]} |B(r)| > \lceil t \rceil^2 \text{ or } \sup_{r \in [t, 2t]} |B(r) - B(t)| > \lceil t \rceil^2\right)$$
$$\leq C \sum_{x \in [-t^2, t^2] \cap \mathbb{Z}} p_x (Y_x \vee Y_{x+1}) + 2e^{-Ct^3}$$
$$\leq CP(\tau_\beta^1(\omega) \geq t, \mathcal{A}_t) Y_0 + 2e^{-Ct^3}.$$

Next, by using $\omega \in \mathcal{B}_1 \cap \mathcal{B}_2 \cap \mathcal{B}_3$, we can replace the logarithm of the probabilities by their $\mathbb{P}$-expectation with the error terms, which yields for $t$ sufficiently large,

$$a_\beta(2t) - (2t)^{\frac{1}{2}+\delta} \leq 2a_\beta(t) + 2t^{\frac{1}{2}+\delta} + \log(1 + 2e^{-t^3 + 2a_\beta(t)}) + C$$
$$\leq 2a_\beta(t) + 2t^{\frac{1}{2}+\delta} + 2e^{-t^3 + 2C(1+t)} + C$$
$$\leq 2a_\beta(t) + Ct^{\frac{1}{2}+\delta},$$

where we have used Lemma 3.1(ii) in the second inequality. This completes the proof of (6.4), and by applying it repeatedly, we obtain for any $k \in \mathbb{N}$,

$$a_\beta(t) \geq \frac{1}{2} a_\beta(2t) - Ct^{\frac{1}{2}+\delta}$$
$$\geq \frac{1}{4} a_\beta(4t) - Ct^{\frac{1}{2}+\delta} 2^{-2(\frac{1}{2}-\delta)} - Ct^{\frac{1}{2}+\delta}$$
$$\geq \cdots$$
$$\geq \left(\frac{1}{2}\right)^k a_\beta(2^k t) - Ct^{\frac{1}{2}+\delta} \sum_{i=0}^{k-1} 2^{-i(\frac{1}{2}-\delta)}.$$

For any $\delta < \frac{1}{2}$, the sum in the last line converges for $k \to \infty$ and we get

$$\frac{a_\beta(t)}{t} + Ct^{-(\frac{1}{2}-\delta)} \geq \lim_{k \to \infty} \frac{a_\beta(2^k t)}{2^k t} = p(\beta). \qquad \square$$

REMARK 6.1. When $d \geq 2$, we have
$$\mathbb{P}(\mathcal{B}_0) = (1 + t^{2d})^{-1}$$
instead of (6.5) and we have to replace (6.6) by
$$\mathbb{P}(\mathcal{B}_1 \cap \mathcal{B}_2 \cap \mathcal{B}_3) \leq Ct^{-2d-1}.$$
This causes no problem since Proposition 2.2 gives us an arbitrarily fast polynomial decay.

PROOF OF PROPOSITION 2.5. Note first that there exists $C > 0$ such that for all $x \in \mathbb{R}$,

(6.7) $$P(B(2) \in dx) \geq CP(B(1) \in dx)e^{Cx^2}.$$

The factor $e^{Cx^2}$ can be regarded as a gain from the 1 extra time. We are going to impose the additional constraint $\{\tau_\infty(\omega) \geq 2\}$ on the left-hand side and show that the additional cost is much smaller than the gain. More precisely, we show that there exists $K(\omega)$ such that for some (deterministic) $c > 0$ and all $x \geq K(\omega)$,

(6.8) $$P^{0,0;2,x}(\tau_\infty(\omega) \geq t) \geq c^{-1} \exp(-c|x|^{\frac{3}{2}}).$$

To see this, denote by $\lambda_k$ the linear function with $\lambda_k(0) = 0$ and $\lambda_k(2) = 5k$, and let $S_k \subseteq \mathbb{R}_+ \times \mathbb{R}$ denote the slanted time-space box

(6.9) $$S_k := \{(s, x) : s \in [0, 2], x \in [\lambda_k(s) - 4, \lambda_k(s) + 4]\}.$$

We write $R_k := |\omega \cap S_k|$ for the number of disasters in $S_k$, and $0 < T_1^{(k)} < \cdots < T_{R_k}^{(k)} < 2$ for the corresponding ordered disaster times. It is convenient to define $T_0^{(k)} := 0$ and $T_{R_k+1}^{(k)} := 2$. As in Section 3, we also consider the interarrival times between disasters:
$$\Delta_i^{(k)} := T_{i+1}^{(k)} - T_i^{(k)} \quad \text{for } i = 0, \ldots, R_k.$$

Note that by our convention $\Delta_0 = T_1^{(k)}$ and $\Delta_{R_k} = 2 - T_{R_k}^{(k)}$. Let us define events
$$\mathcal{E}_k := \{R_k \leq C \log |k|\},$$
$$\mathcal{F}_k := \left\{\min_{i=0,\ldots,R_k} \Delta_i^{(k)} > k^{-\frac{5}{4}}\right\}.$$

Since $R_k$ is Poisson distributed with parameter 8 which has an exponentially decaying tail, we can find $C > 0$ such that
$$\sum_{k \in \mathbb{Z}} \mathbb{P}(\mathcal{E}_k^c) < \infty.$$

Thus we have $\mathbb{P}(\mathcal{E}) = 1$ for $\mathcal{E} := \{\mathcal{E}_k \text{ for all but finitely many } k\}$ by the Borel–Cantelli lemma. Next, note that $\mathcal{F}_k^c$ is nothing but the event that the Poisson process with rate 6 on $[0, 2]$ has a point in the $k^{-5/4}$ neighborhood of the boundary or has two points within distance $k^{-5/4}$. It is easy to see that such probability decays like $\mathbb{P}(\mathcal{F}_k^c) \leq ck^{-5/4}$.

Setting $\mathcal{F} := \{\mathcal{F}_k \text{ for all but finitely many } k\}$ and using the Borel–Cantelli lemma again, we find that $\mathbb{P}(\mathcal{E} \cap \mathcal{F}) = 1$. Now for $\omega \in \mathcal{E} \cap \mathcal{F}$, we find $K(\omega) \geq 2$ such that for all $|k| \geq K(\omega)$, we have $R_k \leq C \log |k|$ and $\min\{\Delta_0^{(k)}, \ldots, \Delta_{R_k}^{(k)}\} > k^{-5/4}$. Observe that every $x \in \mathbb{R}$ is contained in $[5k(x) - 3, 5k(x) + 3]$ for some $k(x) \in \mathbb{Z}$, and in particular $(2, x) \in S_{k(x)}$. Then for all $x$ with $|k(x)| \geq K(\omega)$, we use the estimates from Lemma 3.3 to get

$$P^{0,0;2,x}(\tau_\infty(\omega) \geq t)$$
$$\geq P^{0,0;2,x}(\tau_\beta(\omega) \geq 2, B(u) \in \lambda_{k(x)}(u) + [-3, 3] \text{ for } u \in [0, 2])$$
(6.10)
$$\geq \exp\left(-c - \sum_{i=0}^{R_k} \frac{c}{\Delta_i^{(k)}}\right)$$
$$\geq \exp(-c - c|k|^{\frac{5}{4}} \log |k|).$$

This completes the proof of (6.8).

For $x \in \mathbb{R}$ with $|k(x)| \leq K(\omega)$, we can still use the second line in (6.10) as a lower bound. Therefore, we conclude that

$$P(B(2) \in dx, \tau_\infty(\omega) \geq 2)$$
(6.11)
$$\geq P(B(1) \in dx)\left(C \inf_{x \in \mathbb{Z}^d} e^{Cx^2 - c|x|^{3/2}}\right.$$
$$\left. \wedge \min_{|k| \leq K(\omega)} \exp\left(-c - \sum_{i=0}^{R_k} \frac{c}{\Delta_i^{(k)}}\right)\right). \quad \square$$

## 7. Proof of main result.

PROOF OF THEOREM 1.3. *Part (i)*: This is a direct consequence of Proposition 2.1 and Theorem 2 in [14].

*Part (ii) for $d = 1$*: The almost sure convergence of $t^{-1} \log P(\tau_\infty^1(\omega) \geq t, \mathcal{A}_t)$ to $p(\infty)$ along $t \in \mathbb{N}$ follows by choosing $r = 2$ in Proposition 2.2 and the Borel–Cantelli lemma. Let us extend this convergence to $t \in \mathbb{R}_+$. Note that the definition of the truncation in (1.8) implies $\mathcal{A}_t = \mathcal{A}_{\lceil t \rceil}$. Therefore, we have

$$P(\tau_\infty^1(\omega) \geq \lceil t \rceil, \mathcal{A}_{\lceil t \rceil}) \leq P(\tau_\infty^1(\omega) \geq t, \mathcal{A}_t) \leq P(\tau_\infty^1(\omega) \geq \lfloor t \rfloor, \mathcal{A}_{\lceil t \rceil}).$$

On the other hand, one can easily deduce from Proposition 2.7 and the Borel–Cantelli lemma that almost surely,

$$\log P(\tau^1_\infty(\omega) \geq \lfloor t \rfloor, \mathcal{A}_{\lceil t \rceil}) = \log P(\tau^1_\infty(\omega_{[\lfloor t \rfloor, \lceil t \rceil]^c}) \geq \lceil t \rceil, \mathcal{A}_{\lceil t \rceil})$$
$$\leq \log P(\tau^1_\infty(\omega) \geq \lceil t \rceil, \mathcal{A}_{\lceil t \rceil}) + t^{1/2}$$

for all sufficiently large $t$. Combining the above two bounds, we find

$$\lim_{t \to \infty} \frac{1}{t} \log P(\tau^1_\infty(\omega) \geq t, \mathcal{A}_t) = p(\infty).$$

Next, we get rid of $\mathcal{A}_t$. Since Lemma 3.1 implies $p(\infty) > -\infty$, and since $P(\mathcal{A}_t^c) \leq e^{-ct^3}$, it follows that almost surely,

$$\lim_{t \to \infty} \frac{1}{t} \log P(\tau^1_\infty(\omega) \geq t, \mathcal{A}_t) = \lim_{t \to \infty} \frac{1}{t} \log P(\tau^1_\infty(\omega) \geq t) = p(\infty).$$

Finally, we replace $\tau^1_\infty$ by $\tau_\infty$. By the definition of $\tau^1_\infty(\omega)$, we have

$$\limsup_{t \to \infty} \frac{1}{t} \log P(\tau_\infty(\omega) \geq t) \leq \lim_{t \to \infty} \frac{1}{t} \log P(\tau^1_\infty(\omega) \geq t) \leq p(\infty).$$

On the other hand, using Proposition 2.5, we find that

$$P(\tau_\infty(\omega) \geq t) = \int_\mathbb{R} P^{2,x}(\tau_\infty(\omega) \geq t - 2) P(B(2) \in \mathrm{d}x, \tau_\infty(\omega) \geq 2)$$

(7.1)
$$\geq A(\omega) \int_\mathbb{R} P^{2,x}(\tau_\infty(\omega) \geq t - 2) P(B(1) \in \mathrm{d}x)$$
$$\geq A(\omega) P(\tau^1_\infty(\theta_{1,0}\omega) \geq t - 1).$$

Since $\lim_{t \to \infty} t^{-1} \log P(\tau^1_\infty(\theta_{1,0}\omega) \geq t - 1) = p(\infty)$ almost surely, we are done.

*Part (iii)*: Fix an arbitrary $\delta > 0$. We can use Proposition 2.4 to find $t_0 > 0$ such that for all $\beta \in [0, \infty]$ and $t \geq t_0$,

$$\left| \frac{a_\beta(t)}{t} - p(\beta) \right| \leq \delta.$$

Since $a_\beta(t)$ is the expectation of a random variable depending only on the disasters in a finite area, it is clear that $\beta \mapsto a_\beta(t)$ is continuous. Therefore, there exists $\beta_0$ such that for all $\beta \geq \beta_0$,

$$|p(\beta) - p(\infty)| \leq 2\delta + \frac{1}{t}|a_\beta(t) - a_\infty(t)| \leq 3\delta.$$

This implies the desired continuity. □

7.1. *Almost sure convergence for $d \geq 2$.* In this section, we prove the almost sure convergence $\lim_{t \to \infty} \frac{1}{t} \log P(\tau_\infty(\omega) \geq t) = p(\infty)$ in dimension $d \geq 2$. As mentioned in Remark 2.6, the only point that requires an extra argument is the proof of

$$\text{(7.2)} \qquad \liminf_{t \to \infty} \frac{1}{t} \log P(\tau_\infty(\omega) \geq t) \geq p(\infty).$$

Note that Proposition 2.5 does not generalize to higher dimensions: For $\mathbf{k} \in \mathbb{Z}^d$, let $L_\mathbf{k}$ denote the last disaster in a multi-dimensional version of the time-space box from (6.9). Then almost surely, there exists a point $\mathbf{k} \in \mathbb{Z}^d$ with $|\mathbf{k}| \leq K$ such that $2 - L_\mathbf{k} < K^{-d+1/2}$ for all sufficiently large $K$. If $x$ is behind the last disaster for such $\mathbf{k}$, then in $d \geq 3$ the second line in (6.10) is smaller than $\exp(-cK^{d-1/2}) = o(\exp(-CK^2))$, and hence cannot be compensated by the factor $e^{Cx^2}$ in (6.7). Obviously, the problem in this argument is that we have too many $\mathbf{k}$'s. We solve this problem in the following two steps:

- first restrict $B(1)$ to an essentially one-dimensional slab,
- then show that the above restriction does not affect the limit.

For $\mathbf{k} = (k_2, \ldots, k_d) \in \mathbb{Z}^{d-1}$, let

$$H_\mathbf{k} := \mathbb{R} \times \left[k_2 - \frac{1}{2}, k_2 + \frac{1}{2}\right) \times \cdots \times \left[k_d - \frac{1}{2}, k_d + \frac{1}{2}\right)$$

and set

$$b_t(\omega, \mathbf{k}) := P(\tau_\infty^1(\omega) \geq t, B(1) \in H_\mathbf{k}, \mathcal{A}_t).$$

An easy extension of Proposition 2.5 shows that there exists some positive and finite random variable $A'(\omega)$ such that for all $x \in H_\mathbf{0}$ and $t \geq 2$,

$$P(\tau_\infty(\omega) \geq 2, B(2) \in \mathrm{d}x, \mathcal{A}_t) \geq A'(\omega) P(B(1) \in \mathrm{d}x, \mathcal{A}_t).$$

Then, by the same argument as in (7.1), we have

$$P(\tau_\infty(\omega) \geq t, \mathcal{A}_t) \geq P(\tau_\infty(\omega) \geq t, B(2) \in H_\mathbf{0}, \mathcal{A}_t)$$
$$\geq A'(\omega) b_t(\theta_{1,0}(\omega), \mathbf{0}).$$

Thus (7.2) follows once we show that $\mathbb{P}$-almost surely,

$$\text{(7.3)} \qquad \liminf_{t \to \infty} \frac{1}{t} \log b_t(\omega, \mathbf{0}) \geq p(\infty).$$

The proof of (7.3) is divided into the following two lemmas, which are analogous to Propositions 2.2 and 2.4.

LEMMA 7.1. *There exists $t_0 > 0$ such that for all $t \geq t_0$,*

$$\text{(7.4)} \qquad \mathbb{P}(|\log b_t(\omega, \mathbf{0}) - \mathbb{E}[\log b_t(\omega, \mathbf{0})]| \geq t^{\frac{3}{4}}) \leq t^{-2d-1}.$$

LEMMA 7.2. *There exists $t_0 > 0$ such that for all $t \geq t_0$,*

$$\mathbb{E}[\log b_t(\omega, \mathbf{0})] \geq \mathbb{E}[\log P(\tau_\infty^1(\omega) \geq t, \mathcal{A}_t)] - t^{\frac{3}{4}}. \tag{7.5}$$

PROOF OF LEMMA 7.1. The proof is almost identical to that of Propositions 2.2. Let us introduce a multidimensional version of the notation used before:

$$J_x^{(1)} := x + \left[-\frac{1}{2}, \frac{1}{2}\right)^d,$$

$$J_x^{(5)}(i) := x + 7i\mathbf{e}_1 + \left[-\frac{5}{2}, \frac{5}{2}\right) \times \left[-\frac{1}{2}, \frac{1}{2}\right)^{d-1},$$

$$M^p(\nu) := \sup_{x,y \in \mathbb{R}^d} \min_{i=0,\ldots,p} \nu(J_x^{(5)}(i) \times J_y^{(5)}(i)),$$

where $\mathbf{e}_1, \ldots, \mathbf{e}_d$ denotes the canonical basis of $\mathbb{R}^d$. With these definitions, Lemmas 3.1 and 4.1 readily extend to $d \geq 2$. Moreover, Lemma 5.1 holds for

$$\hat{\nu}_{\omega,\infty}^{r^-,r^+,t}(d(x,y)) := P((B(r^-), B(r^+)) \in d(x,y) \mid \tau_\infty^1(\omega_{[r^-,r^+]^c}) \geq t,$$
$$\mathcal{A}_t, B(1) \in H_\mathbf{0})$$

in place of $\nu_{\omega,\infty}^{r^-,r^+,t}$. Given these ingredients, we can follow the same argument to prove Proposition 2.2.

Let us explain how to verify Lemma 5.1 for $\hat{\nu}$. Since (5.3) holds with $\nu$ replaced by $\hat{\nu}$, the proof of Lemma 5.1 works without change in the case $r^- \geq 2$. The case $r^- < 2$ requires some care because we need to sprinkle the mass on the time interval $[0, 1]$ under the additional constraint $\{B(1) \in H_\mathbf{0}\}$. We define $r_i^+$ as in Section 5 and choose $j_i^+$ such that

$$\hat{\nu}_{\omega,\infty}^{1,r_i^+,t}(H_\mathbf{0} \times J_{j_i^+}^{(1)}) \geq Ct^{-d}. \tag{7.6}$$

Then define $\lambda_i(u) = (\lambda_i^1(u), \lambda^2(u))$ with $\lambda_i^1 : \mathbb{R}_+ \to \mathbb{R}$ and $\lambda_i^2 : \mathbb{R}_+ \to \mathbb{R}^{d-1}$ such that

- $\lambda_i(r_i^+) = j_i^+$,
- $\lambda_i^1$ is linear ($\lambda_i^1(0) = 0$),
- $\lambda_i^2$ is piecewise affine linear with $\lambda_i^2(0) = \lambda_i^2(1) = 0$.

Using this definition, we can replace $\nu$ by $\hat{\nu}$ in (5.9). Observe that, unlike in the one-dimensional case, $S_i^+$ is not a slanted time-space box in the last $d-1$ coordinates. This is in order to ensure $\{B(1) \in H_\mathbf{0}\}$. As a consequence, we have to consider the Brownian bridge conditioned on $\{B(1) \in H_\mathbf{0}\}$ in (5.8). But this does not impose any additional cost since we have the same conditioning in the definition of $\hat{\nu}$. For the coordinates in time and $\mathbf{e}_1$-direction, we can apply an affine transformation and we get (5.9) for $\hat{\nu}$. □

PROOF OF LEMMA 7.2. We argue in a similar way to the proof of Proposition 2.4. Let us introduce events

$$\mathcal{B}_0(t) := \left\{ b_t(\omega, \mathbf{0}) \geq \max_{\mathbf{k} \in \{-t^2, \ldots, t^2\}^{d-1}} b_t(\omega, \mathbf{k}) - e^{-Ct^3} \right\},$$

$$\mathcal{B}_1(t) := \left\{ \left| \log b_t(\omega, \mathbf{0}) - \mathbb{E}[\log b_t(\omega, \mathbf{0})] \right| \leq t^{\frac{3}{4}} \right\},$$

$$\mathcal{B}_2(t) := \left\{ \left| \log P(\tau_\infty^1(\omega) \geq t, \mathcal{A}_t) - \mathbb{E}[\log P(\tau_\infty^1(\omega) \geq t, \mathcal{A}_t)] \right| \leq t^{\frac{3}{4}} \right\}.$$

Note that from here on out $t$ should be replaced by $\lceil t \rceil$, which we omit to ease the notation. Proposition 2.2 and (7.4) yield that for all $t$ large enough

(7.7) $$\mathbb{P}(\mathcal{B}_1(t)^c \cup \mathcal{B}_2(t)^c) \leq 2t^{-2d-1}.$$

Moreover, we claim that

(7.8) $$\mathbb{P}(\mathcal{B}_0(t)) \geq (1 + 2t^2)^{-(d-1)}.$$

Postponing the claim for the moment, note that from (7.7) and (7.8), we get that $\mathcal{B}_0(t) \cap \mathcal{B}_1(t) \cap \mathcal{B}_2(t)$ has a positive probability for all $t$ large enough. In particular, the intersection is not empty and we can choose $\omega \in \mathcal{B}_0(t) \cap \mathcal{B}_1(t) \cap \mathcal{B}_2(t)$. For such an $\omega$, we have

$$\mathbb{E}[\log P(\tau_\infty^1(\omega) \geq t, \mathcal{A}_t)]$$
$$\leq \log P(\tau_\infty^1(\omega) \geq t, \mathcal{A}_t) + t^{\frac{3}{4}}$$
$$= \log \left( \sum_{\mathbf{k} = \{-t^2, \ldots, t^2\}^{d-1}} b_t(\omega, \mathbf{k}) \right) + t^{\frac{3}{4}}$$
$$\leq \log \left( (2t^2 + 1)^d \left( e^{\log b_t(\omega, \mathbf{0})} + e^{-Ct^3} \right) \right) + t^{\frac{3}{4}}$$
$$\leq \log \left( (2t^2 + 1)^d \left( e^{\mathbb{E}[\log b_t(\omega, \mathbf{0})] + t^{\frac{3}{4}}} + e^{-Ct^3} \right) \right) + t^{\frac{3}{4}},$$

where the first and the last inequality follow from $\omega \in \mathcal{B}_1(t) \cap \mathcal{B}_2(t)$, and the second inequality follows from $\omega \in \mathcal{B}_0(t)$. From part (ii) of Lemma 3.1, we see that $\mathbb{E}[\log b_t(\omega, \mathbf{0})]$ decays linearly, which completes the proof of (7.5).

It remains to show (7.8). This is intuitively obvious since $b_t(\omega, \mathbf{0})$ should have the highest chance to be the maximum as it imposes least constraint on $[0, 1]$, and there are $(1 + 2t^2)^{(d-1)}$ many candidates. To make this argument rigorous, it is better to drop the truncation $\mathcal{A}_t$ and work with

$$c_t(\omega, \mathbf{k}) := P(\tau_\infty^1(\omega) \geq t, B(1) \in H_\mathbf{k}).$$

For $\mathbf{k} \in \mathbb{Z}^{d-1}$, let $\omega_\mathbf{k} := \theta_{0,(0,\mathbf{k})}(\omega)$ be obtained by shifting $\omega$ by $(0, \mathbf{k}) \in \mathbb{Z}^d$ in space, and let $\mathbf{K} = \mathbf{K}(\omega)$ be the random index such that

$$c_t(\omega_\mathbf{K}, 0) = \max_{\mathbf{k} = \{-t^2, \ldots, t^2\}^{d-1}} c_t(\omega_\mathbf{k}, 0).$$

Then by this definition, for every $\mathbf{k} \in \{-t^2, \ldots, t^2\}^{d-1}$,

$$c_t(\omega_{\mathbf{K}}, \mathbf{k}) = \int_{H_0} P^{1,z+(0,\mathbf{k})}(\tau^1_\infty(\omega_{\mathbf{K}}) \geq t) P(B(1) \in (0, \mathbf{k}) + dz)$$

$$\leq \int_{H_0} P^{1,z+(0,\mathbf{k})}(\tau^1_\infty(\omega_{\mathbf{K}}) \geq t) P(B(1) \in dz)$$

$$= c_t(\omega_{(0,-\mathbf{k})}, 0)$$

$$\leq c_t(\omega_{\mathbf{K}}, 0).$$

Here, we use $P^{s,x}$ for the law of Brownian motion started at time $s$ with initial distribution $\delta_x$. Together with $P(\mathcal{A}^c_t) \leq e^{-Ct^3}$, we get

$$b_t(\omega_{\mathbf{K}}, 0) \geq c_t(\omega_{\mathbf{K}}, 0) - e^{-Ct^3}$$

$$= \max_{\mathbf{k}=\{-t^2, \ldots, t^2\}^{d-1}} c_t(\omega_{\mathbf{K}}, \mathbf{k}) - e^{-Ct^3}$$

$$\geq \max_{\mathbf{k}=\{-t^2, \ldots, t^2\}^{d-1}} b_t(\omega_{\mathbf{K}}, \mathbf{k}) - e^{-Ct^3}.$$

Now let $\mathbf{L}$ be independent of $\omega$ and uniformly distributed on $\{-t^2, \ldots, t^2\}^{d-1}$, and set $\widetilde{\omega} := \omega_{\mathbf{L}}$. Since $\widetilde{\omega}$ has the same distribution as $\omega$, we have

$$\mathbb{P}(\mathcal{B}_0(t)) = \mathbb{P}\left(b_t(\widetilde{\omega}, 0) \geq \max_{\mathbf{k}=\{-t^2, \ldots, t^2\}^{d-1}} b_t(\widetilde{\omega}, \mathbf{k}) - e^{-Ct^3}\right)$$

$$\geq \mathbb{P}(\mathbf{L} = \mathbf{K}(\omega))$$

$$= (1 + 2t^2)^{-(d-1)}$$

and we are done. □

**Acknowledgments.** Ryoki Fukushima thanks Francis Comets and Nobuo Yoshida for useful discussions on this problem and conveying him the idea of Proposition 1.2. The authors are grateful to Matthias Birkner for pointing out an error in the proof of Lemma 4.1.

RESEARCH INSTITUTE FOR MATHEMATICAL SCIENCES
KYOTO UNIVERSITY
KYOTO
JAPAN
E-MAIL: ryoki@kurims.kyoto-u.ac.jp

DEPARTMENT OF MATHEMATICS
TECHNICAL UNIVERSITY OF MUNICH
MUNICH
GERMANY
E-MAIL: junk@tum.de